\documentclass{article}
\usepackage[utf8]{inputenc}
\usepackage{amsmath, amssymb}
\usepackage{fullpage}
\usepackage{booktabs}
\usepackage{todonotes}
\usepackage{isomath}
\usepackage{bm}
\usepackage{upgreek}
\usepackage{makecell} 
\usepackage{caption, subcaption}
\usepackage{pgfplots}
\pgfplotsset{compat=1.17}

\DeclareMathOperator{\curl}{\mathbf{curl}}
\DeclareMathOperator{\dive}{div}
\DeclareMathOperator{\rot}{rot}
\DeclareMathOperator{\trace}{\mathbf{tr}}

\newcommand{\R}{\mathbb{R}}
\newcommand{\vvh}{\mathbf{v}_h}
\newcommand{\wwh}{\mathbf{w}_h}
\newcommand{\vvhf}{\mathbf{v}_h^F}
\newcommand{\xx}{\mathbf{x}}
\newcommand{\pp}{\mathbf{p}}
\newcommand{\nn}{\mathbf{n}}
\renewcommand{\tt}{\mathbf{t}}
\renewcommand{\P}{\mathbb{P}}
\renewcommand*{\vec}{\vectorsym}

\newcommand{\mat}{\matrixsym}
\newcommand{\dP}{~\text{d}P}
\newcommand{\dF}{~\text{d}F}

\newcommand{\PiE}{{\bm\Pi}^0_{\mathbf{E}}}
\newcommand{\PiB}{{\bm\Pi}^0_{\mathbf{B}}}
\newcommand{\faceProdP}[2]{\left[#1,\,#2\right]_{\text{face},P}}
\newcommand{\edgeProdP}[2]{\left[#1,\,#2\right]_{\text{edge},P}}
\newcommand{\faceProdOm}[2]{\left[#1,\,#2\right]_{\text{face},\Omega}}
\newcommand{\edgeProdOm}[2]{\left[#1,\,#2\right]_{\text{edge},\Omega}}
\newcommand{\ehat}{\hat{\epsilon}}
\newcommand{\shat}{\hat{\sigma}}
\newcommand{\mhat}{\hat{\mu}}

\newcommand{\VEh}{\vec{V}_{h}^{\text{edge}}}
\newcommand{\VBh}{\vec{V}_{h}^{\text{face}}}

\newcommand{\its}[1]{\,\,\,(#1)} 
\title{Parallel block preconditioners for virtual element discretizations of the time-dependent Maxwell equations}
\author{Nicolás A. Barnafi \thanks{Department of Mathematics, Universit\`a degli studi di Pavia, Italy. \texttt{nicolas.barnafi@unipv.it}}\and Franco Dassi \thanks{Dipartimento di Matematica e Applicazioni, Università degli Studi di Milano-Bicocca, Italy.
\texttt{franco.dassi@unimib.it}} \and Simone Scacchi \thanks{Department of Mathematics, Universit\`a degli studi di Milano, Italy. \texttt{simone.scacchi@unimi.it}}}
\date{}
\begin{document}
\maketitle
\begin{abstract}
The focus of this study is the construction and numerical validation of parallel block
preconditioners for low order virtual element discretizations of the three-dimensional Maxwell equations. The virtual element method (VEM) is a recent technology for the numerical approximation of partial differential equations (PDEs), that generalizes finite elements to polytopal computational grids. So far, VEM has been extended to several problems described by PDEs, and recently also to the time-dependent Maxwell equations. When the time discretization of PDEs is performed implicitly, at each time-step a large-scale and ill-conditioned linear system must be solved, that, in case of Maxwell equations, is particularly challenging, because of the presence of both H(div) and H(curl) discretization spaces. We propose here a parallel preconditioner, that exploits the Schur complement block factorization of the linear system matrix and consists of a Jacobi preconditioner for the H(div) block and an auxiliary space preconditioner for the H(curl) block. Several parallel numerical tests have been perfomed to study the robustness of the solver with respect to mesh refinement, shape of polyhedral elements, time step size and the VEM stabilization parameter.
\end{abstract}
\section{Introduction}

Finite element methods (FEM) have been widely used for the numerical approximation of Maxwell equations, see e.g. \cite{nedelec.1980,bermudez.2014,jin.2014,monk.2003}, and the references therein. These theoretical studies have allowed the development of robust solvers, that have been successfully employed for numerical simulations of complex physical phenomena governed by Maxwell equations, such as the design of microwave devices \cite{coccioli.1996}, electromagnetic scattering \cite{khebir.1993}, and antennas \cite{greenwood.1999}.

In many applicative areas of electromagnetism, the computational domain presents complex geometric features, that might require the flexibility of general polytopal grids for an accurate representation or would otherwise requires a large number of simpler elements to be described (tetrahedra and/or hexahedra). The virtual element method (VEM), introduced in \cite{beirao.2013}, consists in a generalization of FEM to the case of polygonal or polyhedral grids constituted by elements of arbitrary shape. So far, VEM has been extended to several problems described by partial differential equations (PDEs), see e.g. \cite{beirao_elast.2013,dassi.2020,beirao_NS.2018,chi.2017,antonietti.2021}. Recently, a low order VEM discretization has been proposed in \cite{fullMax3d} for time-dependent Maxwell equations.

When the time discretization of PDEs is performed implicitly, at each time-step a large-scale and ill-conditioned linear system must be solved. In the case of Maxwell equations, this linear system is particularly challenging, because of the presence of both H(div) and H(curl) discretization spaces. A scalable effective linear solver for FEM discretizations of Maxwell equations has been proposed in \cite{phillips2018scalable}.

In the last five years, researchers have started developing ad-hoc preconditioners for VEM discretizations of different PDE problems: see \cite{antoniettiMasV.2018} for multigrid, \cite{calvo.2019} for Additive Schwarz and \cite{bertoluzzaPP.2020} for BDDC and FETI-DP preconditioners for scalar elliptic equations; see \cite{dassiS.2020b} for parallel block algebraic multigrid (AMG) preconditioners for different saddle point problems; see \cite{dassi.2022} and \cite{bevilacqua.2022} for BDDC preconditioners for mixed formulations of elliptic equations and the Stokes system, respectively. To the best of our knowledge, there are no references yet on effective linear solvers applied to VEM discretizations of Maxwell equations. 

The aim and novelty of the present study is the construction and numerical validation of a parallel block preconditioner for a time-implicit low order VEM discretization of the Maxwell equations in three dimensions. The proposed preconditioner exploits the Schur complement block factorization of the linear system matrix and consists of a Jacobi preconditioner for the H(div) block and an auxiliary space preconditioner for the H(curl) block. Several parallel numerical tests have been perfomed to study the robustness of the solver with respect to mesh refinement, shape of polyhedral elements, time step size and the VEM stabilization parameter. A comparison with a tetrahedral FEM solver is also provided.

\paragraph{Structure of the paper.}
The rest of the manuscript is organized as follows: in Section \ref{sec:cont-model} we introduce Maxwell's equations together with their block structure. In Section \ref{sec:discretization} we show the space and time discretization of the model. We show in detail how the blocks are constructed using VEM. In Section \ref{sec:shur} we show the discrete form of the Schur complement, and how we leverage it to construct a robust preconditioner. The numerical results obtained using such a preconditioner are shown in Section \ref{sec:numerical}. The conclusions of our work are finally given in Section \ref{sec:conclusions}.

\paragraph{Notations.} Throughout this work, we will use the classical Sobolev spaces $H^s(\Omega)$ for $s\in\mathbb N$ and $\Omega\subset \mathbb R^3$ open, bounded and with Lipschitz boundary. 
We endow this space with the standard inner-product $(\cdot, \cdot)_{s,\Omega}$ and 
consider its induced norm $\|\cdot\|=\sqrt{(\cdot, \cdot)_{s,\Omega}}$. 
The particular case of $s=0$ is denoted $L^2(\Omega)$. 
For the definition of the boundary conditions we denote with $H^{1/2}(\partial\Omega)$ the  space of $H^1(\Omega)$ functions restricted to $\partial\Omega$, 
and also consider the space $H^{-1/2}(\partial\Omega)=\left(H^{1/2}(\partial\Omega)\right)'$ its functional dual. In general, for a Banach space X, we will denote its duality product as $\langle\cdot, \cdot\rangle_X$. 
Given a vector function $\vec v=(v_1, v_2, v_3):\Omega\to \R^3$, 
we consider the divergence and curl operators, defined as
$$
\dive \vec v = \partial_x v_1 + \partial_y v_2 + \partial_z v_3, \quad \curl \vec v=(\partial_y v_3 - \partial_z v_2, \partial_z v_1 - \partial_x v_3, \partial_x v_2 - \partial_y v_1), 
$$
which induce the following Sobolev spaces:
\begin{align*}
   H^s( \curl, \Omega) &= \left\lbrace \vec v\in H^s: \curl \vec v \in H^s(\Omega) \right\rbrace, \\
   H^s( \dive, \Omega) &= \left\lbrace \vec v\in H^s: \dive \vec v \in H^s(\Omega) \right\rbrace. 
\end{align*}
If $s=0$, we omit the $s$. 
There exist two trace operators given by $\trace_{\curl}:H(\curl,\Omega)\to \left[H^{-1/2}(\partial\Omega)\right]^2$ and 
$\trace_{\dive}:H(\dive,\Omega)\to H^{-1/2}(\partial\Omega)$ such that, 
if $\vec n$ is the outwards pointing normal vector of $\Omega$, 
$\trace_{\curl} \vec w= \vec w\times \vec n$ and $\trace_{\dive} \vec z = \vec z \cdot \vec n$ for all $\vec w\in H(\curl, \Omega)$ 
and $\vec z \in H(\dive, \Omega)$. With them, the homogeneous boundary condition spaces can be defined as 
\begin{align*}
    H_0(\curl, \Omega) &= \left\lbrace \vec v\in H(\curl, \Omega): \trace_{\curl} \vec v = \vec 0 \right\rbrace, \\
    H_0(\dive, \Omega) &=\left\lbrace \vec v\in H(\dive, \Omega): \trace_{\dive} \vec v = 0 \right\rbrace.
\end{align*}
Moreover, given a generic domain $\mathcal{O}\in\R^d$ ($d=1,2\:$ or 3), 
we refer to the polynomial spaces of degree $k$ defined on $\mathcal{O}$ as $\P_k(\mathcal{O})$.


\paragraph{Notations for meshes.} 
Given a domain $\Omega$ we will denote by $\Omega_h$ a generic partition of $\Omega$ composed by polyhedrons 
whose maximum diameter is $h$.
We refer to a polyhedron of $\Omega_h$ as $P$, 
$\partial P$ is the set of faces that defines the boundary of $P$ and 
$\upzeta P$ is its skeleton, i.e., the set of edges defined by $P$.
Moreover, each face will be denoted by $F$ and 
$\nn_F$ is the unit normal to $F$,
a generic edge is written as $e$ and i
ts direction is denoted by $\tt_e$.
Independently on the geometrical entity we are considering, polyhedron, face, or edge,
we refer to its barycenter via $\xx_\star$.
More specifically, $\xx_P, \xx_F$ and $\xx_e$ refer to the barycenter
of the polyhedron $P$, the face $F$ and the edge $e$, respectively.

\section{The continuous model}\label{sec:cont-model}
We consider the two field formulation of Maxwell's equations, which can be stated as follows: Given the initial data $\vec E^0$ and $\vec B^0$ such that $\dive \vec B^0=0$, find the electric field $\vec E$ and the magnetic induction field $\vec B$ such that   
\begin{equation}\label{eq:maxwell-strong}
\begin{aligned}
   \varepsilon\vec E_t + \sigma \vec E - \curl\left(\mu^{-1}\vec B\right) &= \vec J &\text{ in } \Omega, \forall t\in (0,T),& \\
   \vec B_t + \curl\left(\vec E\right) &= \vec 0 &\text{ in } \Omega, \forall t\in (0,T),& \\
   \vec E(0) &= \vec E^0 &\text{ in } \Omega,& \\
   \vec B(0) &= \vec B^0 &\text{ in } \Omega,& \\
   \trace_{\curl} \vec E &= \vec 0 &\forall t\in (0,T),& \\
   \trace_{\dive} \vec B &= 0 &\forall t\in (0,T),& \\
\end{aligned}
\end{equation}
where the subscript $t$ means the first time derivative, i.e. $\psi_t = \frac{\partial\psi}{\partial t}$ for all functions $\psi$ sufficiently regular. The function $\vec J$ represents the electric current density that is being externally applied, and the parameters are the electric permitivity $\varepsilon$,  the electric conductivity $\sigma$ and the magnetic permeability $\mu$. 

Considering adequate test functions and integrating by parts, the weak formulation of problem \eqref{eq:maxwell-strong} can be obtained: Given the initial data $\vec E^0$ and $\vec B^0$ such that $\dive \vec B^0=0$, find the electric field $\vec E \in H(\curl, \Omega)$ and the magnetic induction field $\vec B\in H(\dive, \Omega)$ such that  

\begin{equation}\label{eq:maxwell-weak}
    \begin{aligned}
      \langle\varepsilon \vec E_t, \vec w\rangle_{H(\curl,\Omega)}  + (\sigma \vec E, \vec w)_{0,\Omega} + (\mu^{-1}\vec B, \curl \vec w)_{0,\Omega} &= \langle\vec J, \vec w\rangle_{H(\curl,\Omega)} &\forall \vec w\in H(\curl, \Omega),& \\
      (\mu^{-1}\vec B_t, \vec \psi)_{0,\Omega} + (\mu^{-1}\psi, \curl \vec E)_{0,\Omega} &= 0 &\forall \vec\psi \in H(\dive, \Omega).&
    \end{aligned}
\end{equation}
We refer to \cite{zhao2004analysis} for specific hypotheses on the data and parameters such that problem \eqref{eq:maxwell-weak} is well-posed. This problem presents a generalized saddle-point structure, and note that the order of the fields is fundamental for having a well-defined Schur complement operator as already hinted in \cite{phillips2018scalable}. Indeed, problem \eqref{eq:maxwell-weak} can be written in operator form as
    \begin{equation}\label{eq:maxwell-block}
        \begin{bmatrix}
            \mathcal A & \mathcal B_1 \\
            \mathcal B_2 & \mathcal C
        \end{bmatrix}
        \begin{bmatrix}
            \vec E \\ \vec B
        \end{bmatrix}
        = 
        \begin{bmatrix}
            \mathcal J \\ \vec 0
        \end{bmatrix}
    \end{equation}
where the corresponding operators are formally given by 
\begin{align*}
    \mathcal A &= \varepsilon\partial_t + \sigma \mathcal I, \\
    \mathcal B_1 &= -\mu^{-1}\curl, \\
    \mathcal B_2 &= \mu^{-1}\curl, \\
    \mathcal C &= \mu^{-1}\partial_t,
\end{align*}
where it holds that $\mathcal B_1^T = \mathcal B_2$. The matrix in \eqref{eq:maxwell-block} can be factorized as
    \begin{equation}\label{eq:schur-factorization}
        \begin{bmatrix}
            \mathcal A & \mathcal B_1 \\
            \mathcal B_2 & \mathcal C
        \end{bmatrix} = 
        \begin{bmatrix}
            \mathcal I & \mat 0 \\
            \mathcal B_2\mathcal A^{-1} & \mathcal I
        \end{bmatrix} 
        \begin{bmatrix}
            \mathcal A & \mat 0\\    
            \mat 0 & \mathcal S_{\mathcal A}
        \end{bmatrix}
        \begin{bmatrix}
           \mathcal I & \mathcal A^{-1} \mathcal B_1 \\
           \mat 0 & \mathcal I 
        \end{bmatrix}
    \end{equation}
where $\mathcal S_{\mathcal A} = \mathcal C - \mathcal B_2\mathcal A^{-1}\mathcal B_1$ is known as the Schur complement. It is possible to define $\mathcal S_{\mathcal C}$ analogously.

\section{Space and time discretization}\label{sec:discretization}

In this section we will give a brief description of the virtual element spaces used 
to have the VEM numerical approximation of the problem defined in Equation~\eqref{eq:maxwell-weak}.
For a more detailed analysis on the interpolation properties of the spaces, 
the underling exact sequence structure and the convergence of fully-discrete scheme,
we refer the reader to~\cite{fullMax3d}.

\subsection{Virtual Element space discretization}\label{sec:vem}

Let $\Omega_h$ be a partition of the three dimensional domain $\Omega$ in generic polyhedral elements. 
To discretize the magnetic and electric fields,
we proceed as the standard virtual element framework:
we define the local spaces and then we glue them together to build the whole functional space.
We describe this procedure for both the electric and magnetic fields in Sections~\ref{sec:ele} and \ref{sec:mag} respectively.

\subsubsection{Electric field -  Edge space}\label{sec:ele}

Consider a generic polyhedron $P\in\Omega_h$,
we define the following space
\begin{align}
\VEh(P):=\bigg\{\vvh\in \big[L^2(P)\big]^3\::\:&\dive(\vvh)=0\,,\quad\curl(\curl(\vvh))\in\big[\P_0(P)\big]^3\,,\nonumber \\[0.3em]
&(\nn_F\times\vvh|_F)\times \nn_F\in\VEh(F),\quad\forall F\in\partial P\,,\nonumber \\[1em]
& \vvh\cdot\tt_e\in\P_0(e),\quad\forall e\in\upzeta P \nonumber \\[1em]
& \int_P \curl(\vvh)\cdot(\xx_P\times \pp_0)\dP=0\quad\forall\pp_0\in[\P_0(P)]^3\:\bigg\}
\label{eqn:edgeLocalP}
\end{align}
where on each face $F\in\partial P$ we have defined the following virtual space 
\begin{align}
\VEh(F):=\bigg\{\vvhf\in \big[L^2(P)\big]^2\::\:&\dive_F(\vvhf)\in\P_0(F)\,,\quad\rot_F(\vvh)\in\P_0(F)\,,\nonumber \\[0.3em]
& \vvhf\cdot\tt_e\in\P_0(e)\quad\forall e\in\partial F,\quad \int_F \vvhf\cdot\xx_F\dF=0\:\bigg\}\,.
\label{eqn:edgeLocalF}
\end{align}
In the definition of $\VEh(F)$ we are considering function $\vvhf$ that are the restriction of $\vvh$ on the face $F$, i.e., 
$\vvhf:=\vvh|_F$.
Moreover, both $\dive_F(\cdot)$ and $\rot_F(\cdot)$ refer to the divergence and rotor differential operators with respect to the two dimensional face coordinate system. 

In such a space a generic function $\vvh$ is uniquely determined by the constant values of $\vvh\cdot\tt_e$ on the skeleton of the polyhedron.
As a consequence the set of degrees of freedom we use is:
\begin{itemize}
\item[\texttt{E1})] the values of $\vvh\cdot\tt_e$ on each edge $e\in\upzeta P$.
\end{itemize}
Moreover, it can be proven that 
constant vector fields belong to $\VEh(P)$, see~\cite{lowest}.
Before defining the global functional space,
we underline that, although a generic function in $\VEh(P)$ is virtual,
it is possible to compute its $L^2$ projection on constant fields.
More specifically, 
it is possible to compute the projection operator $\PiE:\VEh(P)\to[\P_0(P)]^3$ as
\begin{equation}
(\PiE\vvh\,,\pp_0)_P = (\vvh\,,\pp_0)_P\qquad\forall \pp_0\in[\P_0(P)]^3\,,
\label{eqn:L2E}
\end{equation}
starting from the degrees of freedom \texttt{E1} and 
the integral properties plugged in both $\VEh(P)$ and $\VEh(F)$ spaces~\cite{lowest}.
Such projection operator is really important. Indeed, it is one of the key ingredient  to define the bi-linear forms of the discrete weak formulation of the Maxwell's equations.

Finally, the global discrete space for the electric field is
$$
\VEh(\Omega):=\left\{\vvh\in H_0(\curl, \Omega)\::\: \vvh|_P\in\VEh(P)\right\}\,,
$$
where the local spaces are joined together in such a way that 
the tangential component of the vector field on each edge of the mesh is continuous. The degrees of freedom of such a space are obviously the union of the local degrees of freedom \texttt{E1}.  We underline that in this paper we are considering only the lowest order case. It is possible to generalise the definition of such a space for the general order case~\cite{maxwellGO}.

\subsubsection{Magnetic field -  Face space}\label{sec:mag}

Now we describe the virtual element space we use to discretize the magnetic vector field.
Consider a polyhedron $P\in \Omega_h$, we define the following local space
\begin{align}
\VBh(P):=\bigg\{\wwh\in \big[L^2(P)\big]^3\::\:&\dive(\wwh)\in\P_0(P)\,,\quad\curl(\wwh)\in\big[\P_0(P)\big]^3\,,\nonumber \\[0.3em]
&(\nn_F\cdot\wwh|_F)\in\P_0(F),\quad\forall F\in\partial P\,,\nonumber \\[1em]
& \int_P \wwh\cdot(\xx_P\times \pp_0)\dP=0\quad\forall\pp_0\in[\P_0(P)]^3\:\bigg\}\,.
\label{eqn:faceLocalP}
\end{align}

Also $\VBh(P)$ contains not only constant vectorial polynomials but other functions that are virtual.
Thus, a generic function $\wwh\in\VBh(P)$ is uniquely determined by
\begin{itemize}
\item[\texttt{B1})] the face moments on each polyhedron's face
$$
\int_F \wwh\cdot\nn_F\dF\,.
$$
\end{itemize}
Also in this case it is possible to compute the projection operator $\PiB:\VBh(P)\to[\P_0(P)]^3$ as
\begin{equation}
(\PiB\wwh\,,\pp_0)_P = (\wwh\,,\pp_0)_P\qquad\forall \pp_0\in[\P_0(P)]^3\,.
\label{eqn:L2B}
\end{equation}
We refer the reader to~\cite{misto3d} to see how $\PiB$ is computed via the degrees of freedom of $\VBh(P)$.
In Section~\ref{sec:discFrom} we will see that 
also this projection operator is essential to define the bi-linear forms involved in assembling the stiffness matrix.  

Finally, as for the space $\VEh(\Omega)$, 
the global functional space is obtained by gluing each local space together
$$
\VBh(\Omega):=\left\{\wwh\in H_0(\dive, \Omega)\::\: \wwh|_P\in\VBh(P)\right\}\,.
$$
However, the resulting vector field now has 
the normal component continuous across each internal face.
We limit ourselves to the lowest order case, and refer the reader to~\cite{misto3d} for a thorough description of the general case.

\subsection{Discrete weak formulation}\label{sec:discFrom}

Before introducing the discrete counterpart of the problem defined in Equation~\eqref{eq:maxwell-weak},
we have to make two important remarks.

The first one is that $\VEh(\Omega)$ and $\VBh(\Omega)$ are part of an exact sequence~\cite{fullMax3d}.
As a consequence, the following identity holds
\begin{equation}
\curl\left(\VEh(\Omega)\right) = \Big\{\wwh\in\VBh(\Omega)\::\: \dive(\wwh)=0\Big\}\,.    
\label{eqn:divIncl}
\end{equation}
Such relation tells us that if we apply the $\curl$ operator to any function in the space $\VEh(\Omega)$,
we get a function in $\VBh(\Omega)$ whose divergence is zero.

The second remark is the definition of the scalar products within $\VEh(\Omega)$ and $\VBh(\Omega)$ spaces.
More specifically, 
given the space $\VEh(\Omega)$, we define the local discrete counterpart of the $L^2$ inner product as
\begin{equation}
\edgeProdP{\vvh^1}{\vvh^2} := \int_P \PiE\vvh^1\cdot\PiE\vvh^2\dP + 
\mathcal{S}_{\text{edge},P}\left(\left(\mathbf{I}-\PiE\right)\vvh^1,\,\left(\mathbf{I}-\PiE\right)\vvh^2\right)\,,
\label{eqn:innerE}
\end{equation}
where $\vvh^1,\,\vvh^2\in\VEh(P)$ and 
$\mathcal{S}_{\text{edge},P}$ is the stabilization proposed in~\cite[Equation 4.8]{lowest}.
In the same way, the scalar product of the space $\VBh(\Omega)$ is defined as 
\begin{equation}
\faceProdP{\wwh^1}{\wwh^2} := \int_P \PiB\wwh^1\cdot\PiB\wwh^2\dP + 
\mathcal{S}_{\text{face},P}\left(\left(\mathbf{I}-\PiB\right)\wwh^1,\,\left(\mathbf{I}-\PiB\right)\wwh^2\right)\,,
\label{eqn:innerB}
\end{equation}
where now  $\wwh^1,\,\wwh^2\in\VBh(P)$ and
$\mathcal{S}_{\text{face},P}$ is the stabilization proposed in~\cite[Equation 4.17]{lowest}.
Notice that these two scalar products differ by the projection operators ($\PiE$ and $\PiB$) and the stabilization terms ($\mathcal{S}_{\text{edge},P}$ and $\mathcal{S}_{\text{face},P}$). 
Furthermore we underline that,
since the space inclusion defined in Equation~\eqref{eqn:divIncl} holds,
the inner product 
$$
\faceProdP{\curl(\vvh^1)}{\wwh^2}\qquad\text{as well as}\qquad\faceProdP{\wwh^1}{\curl(\vvh^2)}
$$
where $\vvh^1,\vvh^2\in\VEh(P)$ and $\wwh^1,\wwh^2\in\VBh(P)$ are well defined.

Then, once each of these two scalar products are defined on a generic polyhedron $P\in\Omega_h$,
they are defined on the whole domain by summing each polyhedron's contribution, i.e.,
$$
\edgeProdOm{\vvh^1}{\vvh^2} := \sum_{P\in\Omega_h}\edgeProdP{\vvh^1}{\vvh^2} 
\qquad\text{and}\qquad
\faceProdOm{\wwh^1}{\wwh^2} := \sum_{P\in\Omega_h}\faceProdP{\wwh^1}{\wwh^2} 
$$

Starting from the previous observations, 
we are now ready to state the discrete weak formulation of the problem at hand.
Consider a uniform partition of the time interval $[0,T]$ in $n$ sub-intervals of length $\tau$
and an implicit Euler scheme to discretize the time derivative.
Then, the time dependent problem becomes:  given the initial distribution of the fields $\vec E_h^0$ and $\vec B_h^0$, find at each time step $m$ the vector fields $\vec E_h^m$ and $\vec B_h^m$ such that  
\begin{equation}\label{eq:maxwell-fully-discrete}
\left\{
\begin{aligned}
\frac{1}{\tau}\faceProdOm{\mhat^{-1}(\vec B_h^m-\vec B_h^{m-1})}{\wwh} + 
\faceProdOm{\curl(\vec E_h^m)}{\mhat^{-1}\wwh}&=0\\
\frac{1}{\tau}\edgeProdOm{\ehat(\vec E_h^m-\vec E_h^{m-1})}{\vvh} + \edgeProdOm{\shat\vec E_h^m}{\vvh}- 
\faceProdOm{\mhat^{-1}\vec B_h^m}{\curl{\vvh}} &= \edgeProdOm{\bm{J}_I^m}{\vvh}
\end{aligned}
\right.\,,
\end{equation}
hold for each $\wwh\in\VBh(\Omega)$ and $\vvh\in\VEh(\Omega)$.
Here we have defined the element-wise constant functions $\mhat,\ehat$ and $\shat$ 
that approximate the data $\mu,\epsilon$ and $\sigma$, respectively, and 
$J_I^m\in\VEh(\Omega)$ is an interpolation of the current density at the time $t^m$. 

We underline that we swap the equations in the discrete formulation of the problem,
compare Equation~\eqref{eq:maxwell-weak} and Equation~\eqref{eq:maxwell-fully-discrete}.
This choice will become clearer in Section~\ref{sec:shur},
when we define the novel preconditioner proposed in this paper to solve the linear system arising from the discretization of Equation~\eqref{eq:maxwell-fully-discrete}.

\section{Schur complement preconditioner}\label{sec:shur}
The point of departure for the development of an efficient block preconditioner for problem defined in Equation~\eqref{eq:maxwell-fully-discrete} is the Schur complement factorization shown in Equation~\eqref{eq:schur-factorization}, which provides an exact inverse.
Proceeding as in the continuous case, system \eqref{eq:maxwell-fully-discrete} can be written as
    \begin{equation}\label{eq:maxwell-block-discrete}
        \begin{bmatrix}
            \mat C & \mat B_2 \\
            \mat B_1 & \mat A
        \end{bmatrix}
        \begin{bmatrix}
            \vec B_h \\ \vec E_h
        \end{bmatrix}
        = 
        \begin{bmatrix}
            \vec 0\\ \vec J
        \end{bmatrix},
    \end{equation}
where the matrices will be given by
\begin{align*}
    \mat A &= (\varepsilon\tau^{-1} + \sigma) \mat I_{E_h}, \\
    \mat B_1 &= -\mu^{-1}\curl_h, \\
    \mat B_2 &= \mu^{-1}\curl_h, \\
    \mat C &= (\mu\tau)^{-1}\mat I_{B_h},
\end{align*}
with $\mat B_1^T = \mat B_2$. More specifically, The inverse of \eqref{eq:maxwell-block-discrete} is given by
$$
\begin{bmatrix}
  \mat I_{B_h} & -\mat C^{-1}\mat B_2 \\ \mat 0 & \mat I_{E_h}  
\end{bmatrix}
\begin{bmatrix}
 \mat C^{-1} & \mat 0 \\ \mat 0 & \mat S_{\mat C}^{-1} 
\end{bmatrix}
\begin{bmatrix}
\mat I_{B_h} & \mat 0 \\ -\mat B_1\mat C^{-1} & \mat I_{E_h} 
\end{bmatrix},
$$
where the difficulty for the development of a preconditioner has been reduced to the development of efficient preconditioners for the operators $\mat C$ and $\mat S_{\mat C}$, the latter being given by 
    $$ \mat S_{\mat C} = \mat A - \mat B_1\mat C^{-1} \mat B_2 = (\varepsilon\tau^{-1}+\sigma)\mat I_{E_h} + \tau\mu^{-1}\curl_h \mat I_{B_h}\curl_h. $$

The last part of the preconditioner $\curl_h \mat I_{B_h}\curl_h$ is built combining the definition of the face inner product Equation~\eqref{eqn:innerB} and commuting diagram property.
For more detail about we refer to Proposition 6 in~\cite{fullMax3d}.

This yields two fundamental points: on one hand, $\mat S_{\mat C}$ belongs to $(\VEh)'$ and so it is well defined only on an edge functional space, which justifies the change of order considered for the problem variables. On the other hand, it is a differential operator that defines a simple bilinear form, and so it allows for a sparse representation for the Schur complement, which is ideal for the development of a preconditioner. In what follows we argue our choices of preconditioners for the matrices $\mat C$ and $\mat S_{\mat C}$.

\paragraph{Preconditioner for $\mat C$.} This matrix is a scaled mass matrix, for which a simple diagonal preconditioner provides good performance in FEM formulations \cite{wathen1987realistic}. It is still not clear whether a diagonal preconditioner is a good preconditioner in VEM, so this will be evaluated as well in our study.
\paragraph{Preconditioner for $\mat S_{\mat C}$.} The preconditioning of the operator $I + \curl\curl$ has been most successful through auxiliary space methods within multigrid methods \cite{kolev2009parallel}. This has been implemented in HYPRE \cite{falgout2002hypre} in the AMS solver (Auxiliary-space Maxwell solver), whose setup requires the coordinates of the degrees of freedom and a discrete gradient given by an orientation of the mesh edges. 

We highlight that Schur complement preconditioners are attractive because the preconditioned problem has at most four distinct eigenvalues \cite{murphy2000note} in exact arithmetic. This is often not observed in practice unless a sparse representation of the Schur complement $\mat S_{\mat C}$ is available, which in fact is our case. 
\section{Numerical tests}\label{sec:numerical}
The scope of this section is twofold. On one hand, we want to explore the parallel performance of state-of-the-art preconditioners for VEM and their dependence on the timestep. On the other hand, we want to present a first comparison of the solution between FEM and VEM, for an increasing level of complexity of the meshes under consideration. For the FEM problem we consider the Firedrake library \cite{rathgeber2016firedrake}, which already supports an efficient setup of the AMS preconditioner for the $I + \curl\curl$ operator, and for the VEM problem we consider 5 meshes, starting from a comparable tetrahedral mesh to more complex ones. 

\begin{figure}
    \centering
    \begin{subfigure}{0.49\textwidth}
        \includegraphics[width=\textwidth]{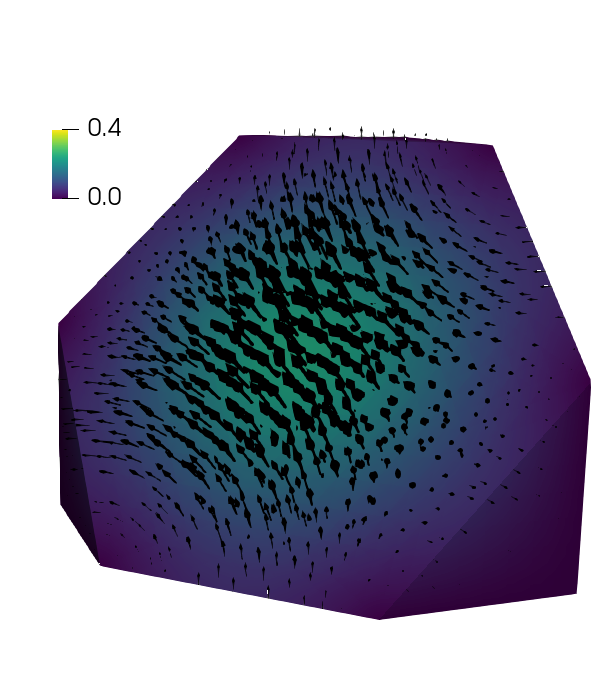}
        \caption{Electric field $E_h$.}
    \end{subfigure}
    \begin{subfigure}{0.49\textwidth}
        \includegraphics[width=\textwidth]{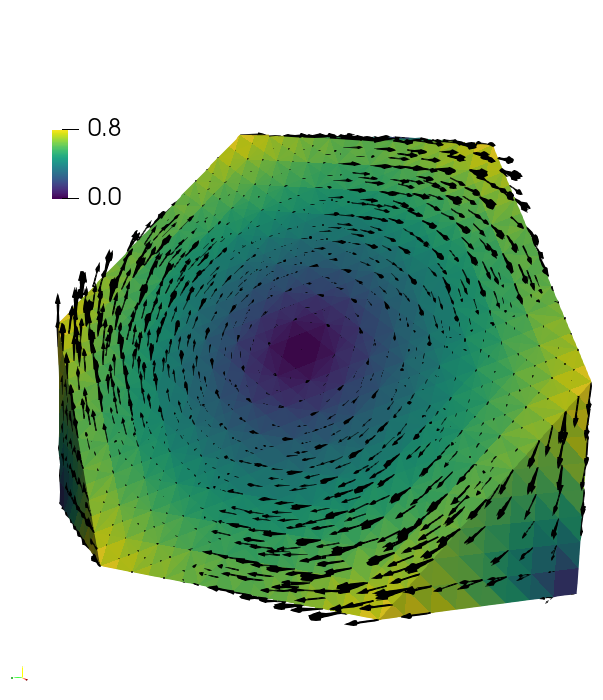}
        \caption{Magnetic field $B_h$.}
    \end{subfigure}
    \caption{Solution of the problem at instant $t=0.5$, where the geometry has been clipped by a plane of normal $(1,1,1)$ passing by the point $(0.5,0.5,0.5)$. In both figures, the color represents the magnitude of the vector being drawn.}
    \label{fig:solution}
\end{figure}

The problem setup is the same in all cases: homogeneous initial and boundary conditions and $\vec J=(1,1,1)$. The solver is given by a GMRES method without restart, the absolute and relative tolerances are given by $10^{-12}$ and $10^{-6}$ respectively, with divergence defined as more than 1000 iterations required for convergence. The block preconditioner was set using the PCFIELDSPLIT preconditioner from PETSc \cite{petsc-user-ref} with a lower-triangular factorization, and the sub-block preconditioners described in Section \ref{sec:shur}. The implementation was done in an {\it in-house} VEM library, and the solution with tetrahedra is shown in Figure \ref{fig:solution}.

\subsection{Meshes under consideration}\label{sec:meshes}
\begin{figure}[!htb]
    \centering
    \begin{tabular}{ccc}
    \includegraphics[width=0.30\textwidth]{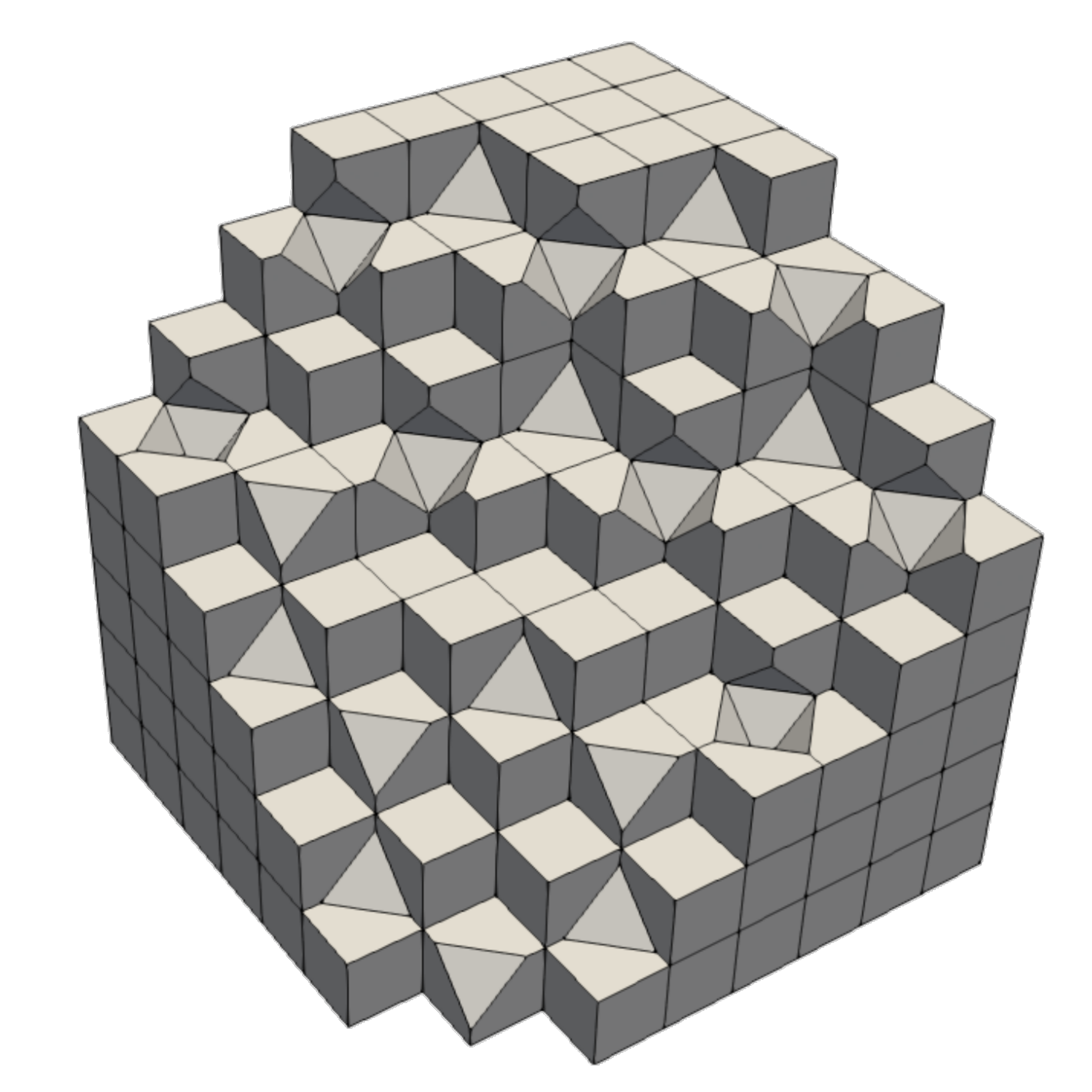}&
    \includegraphics[width=0.30\textwidth]{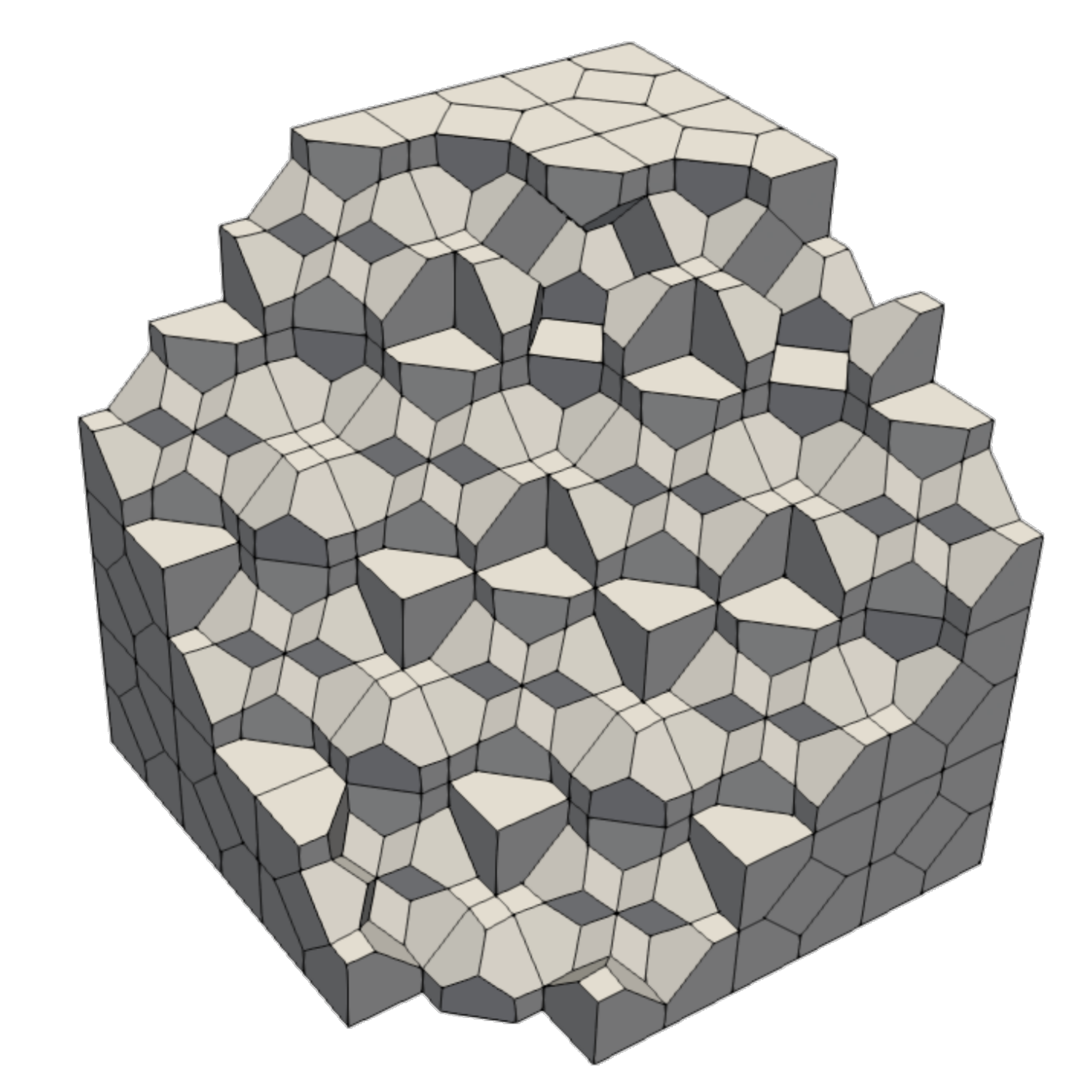}&
    \includegraphics[width=0.30\textwidth]{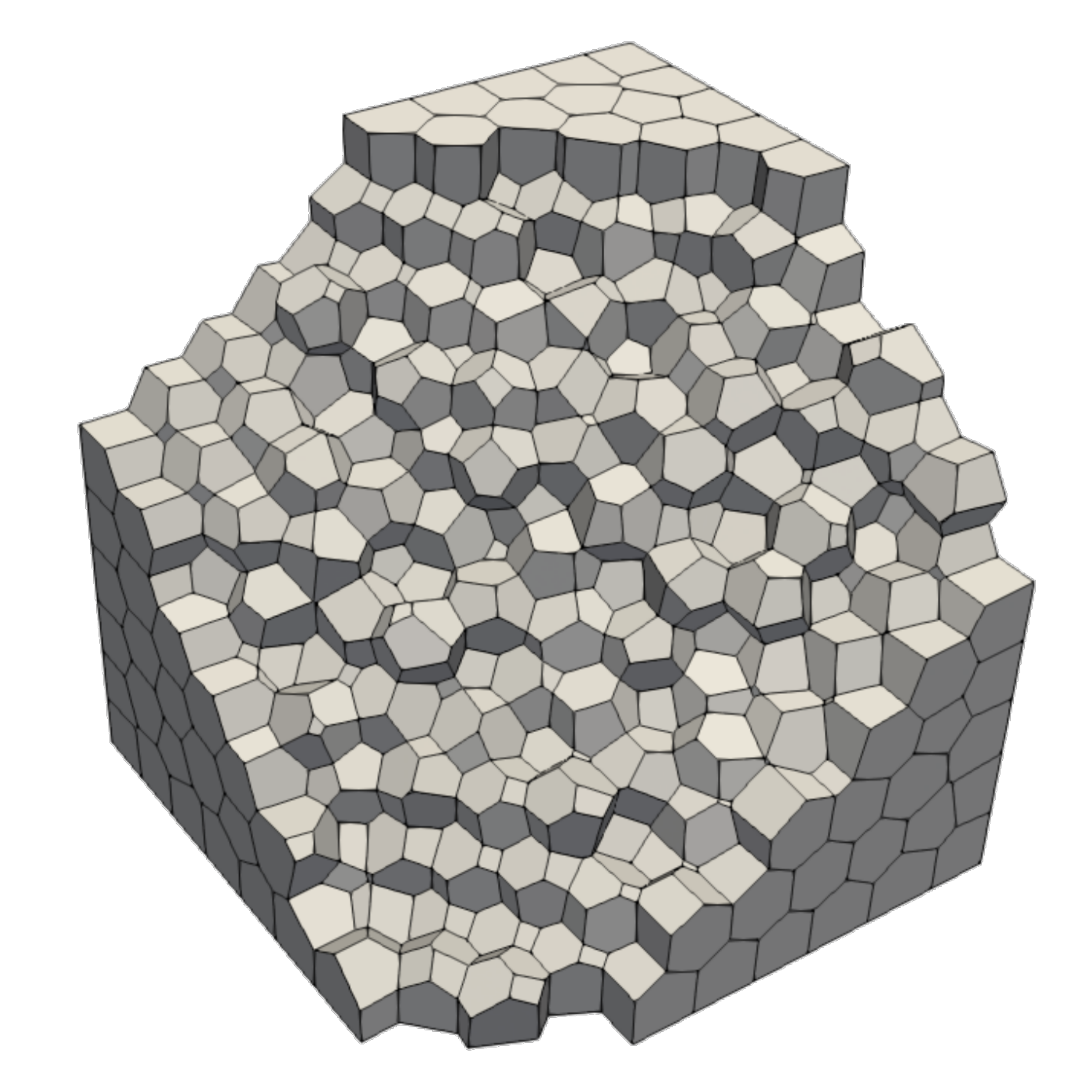}\\
    Octahedral &Nonahedral &Irregular polyhedral
    \end{tabular}
    \caption{Types of meshes used in the numerical experiments}
    \label{fig:meshes}
\end{figure}

Since the virtual element method can deal with polygons of arbitrary shape,
we consider different types of meshes 
to show the performance of the proposed preconditioners.
More specifically, we build the following type of meshes:
\begin{itemize}
\item \textbf{Tetrahedral}: a standard tetrahedral mesh built via the Delaunay criterion~\cite{tetgen};
\item \textbf{Hexahedral structured}: a standard mesh composed by cubes; 
\item \textbf{Octahedral}: a mesh composed by a stencil of eight elements repeated to cover the whole domain;
\item \textbf{Nonahedral}: a mesh composed by another stencil of nine elements repeated to cover the whole~domain;
\item \textbf{Irregular polyhedral}: centroidal voronoi tesselation of the domain optimized via a Lloyd algoritm~\cite{voroPlusPlus}.
\end{itemize}
Since the first two types of meshes are standard,
we show only an examples of Octahedral, Nonahedral and Irregular polyhedral types in Figure~\ref{fig:meshes}.

The meshes taken into account have an increasing complexity.
Tetrahedral and Hexahedral structured are standard meshes composed by regularly shaped elements,
i.e., high quality tetrahedrons and standard cubes, respectively.
Octahedral and Nonhedral meshes are composed by polygons but,
since they are composed by the same stencil,
polyhedrons' shape is regular, they do not present small edges and faces.
Finally, Irregular polyhedral is the most challenging mesh type among the ones taken into account.
Indeed, polyhedrons are different one to each other and additionally have small edges and faces.
This last mesh is expected to be the most difficult to solve with an iterative method since it presents the worst conditioning.
To have a numerical evidence about this fact we compute via Matlab the conditioning of a small matrix (approximately 1500 DoFs) 
arising from the virtual element discretization of the problem at hand in Table~\ref{tab:cond}. This is also relevant to understand the elevated number of linear iterations incurred in our numerical tests in the tetrahedral mesh, which at first could seem like the simplest one.
\begin{table}[!htb]
    \centering
    \begin{tabular}{c|c|c|c|c}
    \toprule
    Tetrahedral &Hexahedral structured &Octahedral &Nonahedral &Irregular polyhedral\\
    \midrule
     $5.53 \cdot 10^2$  & $1.67\cdot 10^1$  &$5.29 \cdot 10^1$  & $7.46\cdot 10^1$ &$1.21 \cdot 10^3$\\
    \bottomrule
    \end{tabular}
    \caption{Conditioning of a small matrix (approximately 1500 DoFs) 
arising from the virtual element discretization of Maxwell's equations varying mesh types.}
    \label{tab:cond}
\end{table}

\subsection{Results: FEM}
In this section we present the performance of an optimal preconditioner for FEM as a point of reference to evaluate the performance of the same scenario with VEM. In other words, we show the performance of the preconditioner described in Section~\ref{sec:shur}, whose optimal scalability has already been studied for FEM in \cite{phillips2018scalable}.

In Table \ref{tab:optimality-fem} we report the performance of the preconditioner for up to more than $6\cdot 10^5$ degrees of freedom for various time steps. As expected, the preconditioner is robust with respect to the problem size, and its performance improves when we reduce the time step. The strong scalability of the preconditioner is reported in Table \ref{tab:scalability-fem}, where it can be observed that there only a small increase in the number of GMRES iterations, and the solution times present an adequate reduction as the number of CPU cores is increased. In general, our results are in agreement with the ones obtained in \cite{phillips2018scalable}.

\begin{table}
    \centering
    \small
        \begin{tabular}{r|r|r|r|r}
            \toprule Dofs & \multicolumn{4}{c}{Solution time (linear iterations)} \\
         \midrule & \makecell[c]{$\Delta t=0.1$} & \makecell[c]{$\Delta t=0.05$} & \makecell[c]{$\Delta t=0.01$} & \makecell[c]{$\Delta t=0.005$}  \\
             \cmidrule(lr){2-5} 218 & $2.80\,10^{-1} \its{14}$  & $2.50\,10^{-1} \its{11}$ & $2.70\,10^{-1} \its 9$ & $2.60\,10^{-1} \its 9$ \\
         1464 & $2.90\,10^{-1} \its{16}$ & $2.80\,10^{-1} \its{14} $ & $2.70\,10^{-1} \its 9$ & $2.60\,10^{-1} \its 9$ \\
         10712 &  $5.10\,10^{-1} \its{18}$ & $4.90\,10^{-1} \its{17}$ & $4.40\,10^{-1} \its{11}$ & $4.40\,10^{-1} \its{10}$ \\
         81712 & $1.85\,10^0 \its{18}$ &  $1.79\,10^0 \its{18}$ & $1.29\,10^0 \its{14}$ & $1.15\,10^0 \its{11}$ \\
         638048 & $1.04\,10^1 \its{18}$ & $1.04\,10^1 \its{18}$ & $7.92\,10^0 \its{17}$ & $7.10\,10^0 \its{14}$ \\ \bottomrule
        \end{tabular}
    \caption{Optimality test, structured tetrahedral mesh with FEM using the optimal preconditioner.}
    \label{tab:optimality-fem}
\end{table} 
\begin{table}
    \centering
    \small
        \begin{tabular}{r|r|r|r|r|r}
            \toprule \multicolumn{6}{c}{Solution time (linear iterations)} \\
         \midrule DoFs & \makecell[c]{1 CPU} & \makecell[c]{2 CPU} & \makecell[c]{4 CPU} & \makecell[c]{8 CPU} & \makecell[c]{16 CPU}  \\
         \cmidrule(lr){1-6} 638048 & $7.33\,10^{0} \its{18}$ & $4.15\,10^0 \its{19}$ & $2.46\,10^0 \its{20}$ & $1.83\,10^0 \its{20}$ & $1.58\,10^0 \its{20}$ \\ \bottomrule
    \end{tabular}
    \caption{Scalability test, FEM.}
    \label{tab:scalability-fem}
\end{table} 


\subsection{Results: VEM} 
In this section we study the performance of the proposed preconditioners on all of the meshes displayed in Section \ref{sec:meshes}. We focus on three aspects of the solvers: (i) their optimality, (ii) their strong scalability, and (iii) their robustnes with respect to the stabilization term detailed in Section \ref{sec:vem}.

\subsubsection{Optimality tests}
In this section we report the robustness of the proposed preconditioner with respect to the problem size. To obtain a better intuition on which block is the one yielding a higher contribution to the deterioration of the preconditioner, as well as validating the efficiency of the sparse representation of the Schur complement, we present the results in four different scenarios: 
\begin{itemize}
    \item {\bf Exact-exact}: both $\mat C$ and $\mat S_{\mat C}$ are solved exactly. 
    \item {\bf Jacobi-exact}: the inverse of $\mat C$ is approximated by the action of a Jacobi preconditioner and $\mat S_{\mat C}$ is solved exactly (with an LU decomposition).
    \item {\bf Exact-AMS}: Matrix $\mat C$ is inverted exactly, whereas the inverse of $\mat S_{\mat C}$ is approximated by the action of the AMS preconditioner.
    \item {\bf Jacobi-AMS}: The inverse of both matrices $\mat C$ and $\mat S_{\mat C}$ and approximated by the application of the previously described preconditioners.
\end{itemize}

\begin{figure}
    \pgfplotsset{every axis plot/.append style={line width=1.5pt}}
    \pgfplotsset{every axis plot/.append style={mark options={draw opacity=0.0, fill opacity=1.0, mark size=1.8pt, mark=square*, }} }
    \centering
    \begin{subfigure}{0.49\textwidth}
        \begin{tikzpicture}
    \begin{loglogaxis}
        [width=\textwidth,
        xlabel=DoFs,
        ylabel=Iterations,
         legend columns=2, 
         legend pos=north west
        ]

        \addplot  coordinates { 
            (349,19) 
            (1404,25) 
            (9156,33) 
            (72422,43) 
            (648444,67) 
            };
        \addplot  coordinates { 
            (349,18) 
            (1404,24) 
            (9156,30) 
            (72422,36) 
            (648444,44) 
            };
         \addplot  coordinates { 
            (349,15) 
            (1404,20) 
            (9156,26) 
            (72422,28) 
            (648444,30) 
            };                               
         \addplot  coordinates { 
            (349,14) 
            (1404,19) 
            (9156,24) 
            (72422,27) 
            (648444,28) 
            };           
        \addplot [mark=none, color=green!60!black] coordinates { 
            (1404,12 * 1.5^2) 
            (9156,12 * 1.5^3) 
            (72422,12 * 1.5^4) 
            (648444,12 * 1.5^5) 
            };
    \legend{0.1, 0.05, 0.01, 0.005, linear}
    \end{loglogaxis}
    \end{tikzpicture}
    
    \caption{Tetrahedral mesh}
    \end{subfigure}
    \begin{subfigure}{0.49\textwidth}
        \begin{tikzpicture}
    \begin{loglogaxis}
        [width=\textwidth,
        xlabel=DoFs,
        ylabel=Iterations,
         legend columns=2, 
         legend pos=north west 
        ]
         
        \addplot  coordinates { 
            (90,2)  
            (540,7)  
            (3672,13)  
            (26928,20)  
            (205920,32)  
            };
        \addplot  coordinates { 
            (90,2)  
            (540,6)  
            (3672,10)  
            (26928,14)  
            (205920,21)  
            };
         \addplot  coordinates { 
            (90,2)  
            (540,4)  
            (3672,7)  
            (26928,9)  
            (205920,11)  
            };                               
         \addplot  coordinates { 
            (90,2)  
            (540,4)  
            (3672,7)  
            (26928,9)  
            (205920,10)  
            };           
            \addplot [mark=none, color=green!60!black] coordinates { 
            (540,5 * 1.5^2) 
            (3672,5 * 1.5^3) 
            (26928,5 * 1.5^4) 
            (205920,5 * 1.5^5) 
            };
    \legend{0.1, 0.05, 0.01, 0.005, linear}
    \end{loglogaxis}
    \end{tikzpicture}
    
    \caption{Hexahedral mesh}
    \end{subfigure}
    
    \begin{subfigure}{0.49\textwidth}
        \begin{tikzpicture}
    \begin{loglogaxis}
        [width=\textwidth,
        xlabel=DoFs,
        ylabel=Iterations,
         legend columns=2, 
         legend pos=north west
        ]
         
        \addplot  coordinates { 
            (110,8)
            (700,12)
            (4952,16)
            (37168,24)
            (237690,40)
            };
        \addplot  coordinates { 
            (110,7)
            (700,11)
            (4952,14)
            (37168,17)
            (237690,25)
            };
         \addplot  coordinates { 
            (110,6)
            (700,9)
            (4952,11)
            (37168,12)
            (237690,13)
            };                               
         \addplot  coordinates { 
            (110,6)
            (700,8)
            (4952,10)
            (37168,11)
            (237690,12)
            };           
        \addplot [mark=none, color=green!60!black] coordinates { 
            (700,6 * 1.5^2) 
            (4952,6 * 1.5^3) 
            (37168,6 * 1.5^4) 
            (237690,6 * 1.5^5) 
            };
    \legend{0.1, 0.05, 0.01, 0.005, linear}
    \end{loglogaxis}
    \end{tikzpicture}
    
    \caption{Octahedral mesh}
    \end{subfigure}
    \begin{subfigure}{0.49\textwidth}
        \begin{tikzpicture}
    \begin{loglogaxis}
        [width=\textwidth,
        xlabel=DoFs,
        ylabel=Iterations,
        legend columns=2,
        legend pos=north west
        ]
         
        \addplot  coordinates { 
            (209,13) 
            (1408,17) 
            (10376,24) 
            (79600,57) 
            };
        \addplot  coordinates { 
            (209,12) 
            (1408,15) 
            (10376,20) 
            (79600,48) 
            };
         \addplot  coordinates { 
            (209,10) 
            (1408,12) 
            (10376,15) 
            (79600,32) 
            };                               
         \addplot  coordinates { 
            (209,9) 
            (1408,11) 
            (10376,13) 
            (79600,28) 
            };           
         \addplot [mark=none, color=green!60!black] coordinates { 
            (1408,13 * 1.5^2) 
            (10376,13 * 1.5^3) 
            (79600,13 * 1.5^4) 
            };
    \legend{0.1, 0.05, 0.01, 0.005, linear}
    \end{loglogaxis}
    \end{tikzpicture}
    
    \caption{Nonahedral mesh}
    \end{subfigure}
    \caption{Optimality test of the proposed preconditioner for time steps $0.1, 0.05, 0.01$ and $0.005$.}
    \label{fig:optimality}
\end{figure}
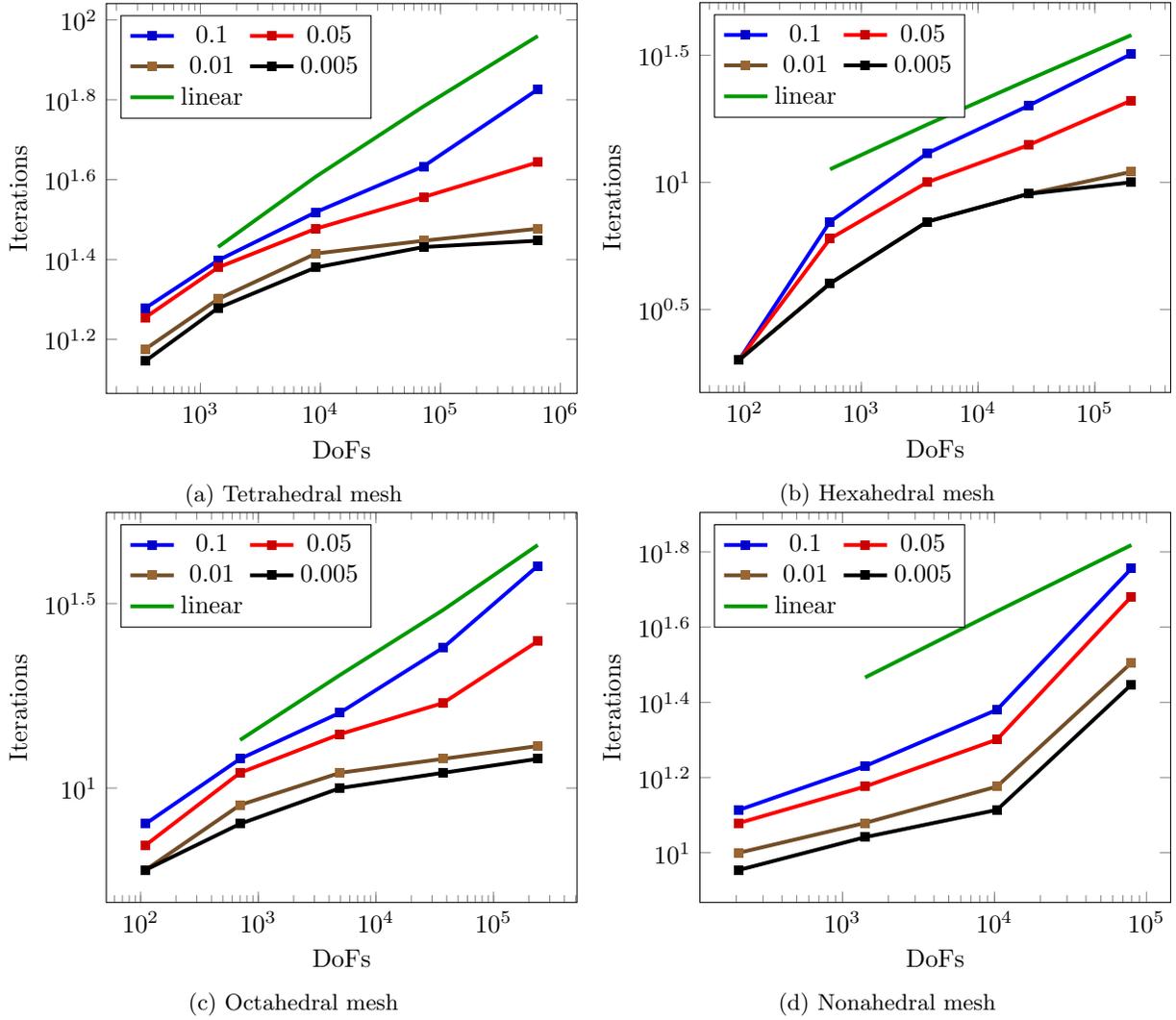

The results of the Jacobi-AMS preconditioner are displayed, for all timesteps considered, in Figure \ref{fig:optimality}. In particular, we can appreciate that for all but the biggest timestep, there is a sublinear growth in the number of linear iterations. This is not true for the nonahedral mesh, where the deterioration in performance is indeed worse that linear. In what follows, we provide further insight into the results obtained for each mesh separately. 

\paragraph{Tetrahedral mesh.} We show the results in Table \ref{tab:optimality-tetgen}. We first note that the Schur complement operator is represented accurately, as GMRES incurrs on only 2 iterations when both inverses are computed exactly. Both preconditioners present a deterioration when using the largest time step (0.1), but AMS is closer to being linear, meaning that it is more sensitive to the time step size. This effect is alleviated when considering a smaller time step, as can be seen by the performance of the combined Jacobi-AMS preconditioner for the second largest time step, which presents an increase of only a couple of iterations for each level of refinement. It is interesting to note that the Jacobi preconditioner itself requires more iterations to converge, but its increase in iterations is milder than the AMS. For example, for $\tau=0.05$, the factor by which the iterations increase between the coarsest and finest mesh is 2.33, whereas the same for the AMS preconditioner is 5.33. This shows that performance is mainly hindered by the AMS preconditioner in the largest time step scenario. All other smaller time steps present a similar number of linear iterations between the Jacobi-AMS and the Jacobi-exact preconditioners.

\paragraph{Hexahedral structured mesh.} The results are shown in Table \ref{tab:optimality-struct}, which are similar to the ones obtained with the tetrahedral mesh. The main difference between this scenario and the previous one is that this one is much easier as shown in Table \ref{tab:cond}, which is further confirmed by the overall reduced number of iterations. The deterioration in performance of both Jacobi and AMS preconditioners is comparable as shown in the subtables (b) and (c) of Table \ref{tab:optimality-struct}, and while their performance is still not quite robust for $\tau= 0.1$, it suffices to consider a smaller time step (0.05) in order to obtain a robust preconditioner. We note that, as in the tetrahedral case, for all but the largest time step, the number of linear iterations is similar between the Jacobi-exact and Jacobi-AMS preconditioners, meaning that AMS is the main driver of deterioration for larger time steps. 

\paragraph{Octahedral mesh.} This is the first case considered with a truly polygonal mesh, and its results are shown in Table \ref{tab:optimality-octa}. This example presents similar behaviour to the previous ones, with the performance being mainly hindered by the AMS preconditioner for the largest time step. Instead, the number of iterations is similar between Jacobi-exact and Jacobi-AMS.

\paragraph{Nonahedral mesh.} The results of this case are shown in Table \ref{tab:optimality-novedri}. This test is interesting as its trends change with the time step. More specifically, it is still true that the performance of the Exact-AMS preconditioner presents deterioration for the largest time step, but this instead does not hold for smaller time steps. Instead, in all cases there is a significant deterioration of the Jacobi preconditioner, as shown by the performance of the Jacobi-exact, on the finest mesh. As before, the Jacobi-AMS preconditioner presents similar number of iterations to the Jacobi-exact preconditioner, but as the latter is not robust in this case, neither is the former. 

\paragraph{Irregular polyhedral mesh.} As expected, this case is the most difficult, with the results being shown in Table \ref{tab:optimality-voro}. The ill-conditioning of the problem is severe for this type of mesh, which can be noted by the increase in GMRES iterations in the finest mesh for the largest time step (4 iterations instead of 2), despite the use of an exact inverse. In contrast to the nonahedral case, the Jacobi preconditioner is able to adequately handle the mass matrix, and instead the AMS preconditioner is not able to capture the physics, as observed by its poor convergence and even worse robustness with respect to the problem size, as shown by its divergence in most cases. This can be mildly alleviated by reducing the time step, but the overall time reduction is not significant and not worth the effort due to the lack of robustness.\\[0.3em] 

\noindent We highlight that in all cases, an adequate performance is shown by the sparse representation of the Schur operator $\mat S_{\mat C}$, and in all cases but the irregular one, a slightly smaller time step results in a preconditioner with an adequate performance for large scale simulations. It is also fundamental to note that in all such cases, there is a difference of up to two orders of magnitude between the exact solver (Exact-exact) and the iterative one (Jacobi-AMS).

\begin{table}
    \centering
    \small
    \begin{subtable}{0.6\textwidth}
        \begin{tabular}{r|r|r|r|r}
            \toprule Dofs & \multicolumn{4}{c}{Solution time (linear iterations)} \\
        \midrule & \makecell[c]{$\Delta t=0.1$} & \makecell[c]{$\Delta t=0.05$} & \makecell[c]{$\Delta t=0.01$} & \makecell[c]{$\Delta t=0.005$} \\  
         \cmidrule(lr){2-5}  349 & $1.04\,10^{-2} \its 2 $   & $7.57\,10^{-3} \its 2$  & $9.52\,10^{-3} \its 2 $ & $7.46\,10^{-3} \its 2 $ \\
         1404 & $9.59\,10^{-3} \its 2$ & $9.81\,10^{-3} \its 2 $ & $1.02\,10^{-2} \its 2$ & $9.66\,10^{-3} \its 2$ \\
         9156 & $7.51\,10^{-2} \its 2$ & $7.92\,10^{-2} \its 2$ & $7.75\,10^{-2} \its 2$ & $7.63\,10^{-2} \its 2$ \\
         72422 & $8.15\,10^0 \its 2$ &  $8.41\,10^0 \its 2$ & $8.40\,10^0\its 2$ & $8.31\,10^0 \its 2$ \\
         648444 & $1.14\,10^3 \its 2$ & $1.16\,10^3 \its 2$ & $1.14\,10^3 \its 2 $ & $1.17\,10^3 \its 2 $ \\ \bottomrule
        \end{tabular}
        \caption{Exact-exact}
    \end{subtable}
    \begin{subtable}{0.6\textwidth}
        \begin{tabular}{r|r|r|r|r}
            \toprule Dofs & \multicolumn{4}{c}{Solution time (linear iterations)} \\
         \midrule & \makecell[c]{$\Delta t=0.1$} & \makecell[c]{$\Delta t=0.05$} & \makecell[c]{$\Delta t=0.01$} & \makecell[c]{$\Delta t=0.005$} \\ 
         \cmidrule(lr){2-5}  349 & $7.07\,10^{-3} \its{19}$  & $7.15\,10^{-3} \its{18}$ & $7.02\,10^{-3} \its{15}$ & $6.77\,10^{-3} \its{15}$ \\
         1404 &  $9.44\,10^{-3} \its{25}$ & $9.72\,10^{-3} \its{24}$ & $9.46\,10^{-3} \its{20}$ & $9.43\,10^{-3} \its{19}$ \\
         9156 & $6.81\,10^{-2} \its{32}$ & $6.85\,10^{-2} \its{29}$ & $6.40\,10^{-2} \its{25}$ & $6.18\,10^{-2} \its{24}$ \\
         72422 &  $6.81\,10^{-2} \its{32}$ &  $6.85\,10^{-2} \its{29}$ & $6.40\,10^{-2} \its{25}$ & $6.18\,10^{-2} \its{24}$ \\
         648444 & $8.47\,10^2 \its{49}$ & $1.16\,10^2 \its{42}$ & $1.14\,10^2 \its{29}$ & $1.17\,10^2 \its{27}$ \\ \bottomrule
        \end{tabular}
        \caption{Jacobi-exact}
    \end{subtable}   
    \begin{subtable}{0.6\textwidth}
        \begin{tabular}{r|r|r|r|r}
            \toprule Dofs & \multicolumn{4}{c}{Solution time (linear iterations)} \\
          \midrule & \makecell[c]{$\Delta t=0.1$} & \makecell[c]{$\Delta t=0.05$} & \makecell[c]{$\Delta t=0.01$} & \makecell[c]{$\Delta t=0.005$} \\ 
         \cmidrule(lr){2-5} 349 & $7.69\,10^{-3} \its{3}$  & $8.00\,10^{-3} \its{3}$ & $7.98\,10^{-3} \its 3$ & $7.70\,10^{-3} \its 3$ \\
         1404 & $1.22\,10^{-2} \its 5$ & $1.24\,10^{-2} \its 4$ & $1.22\,10^{-2} \its 4$ & $1.27\,10^{-2} \its 4$ \\
         9156 & $7.07\,10^{-2} \its 8$ & $6.88\,10^{-2} \its 7$ & $6.88\,10^{-2} \its 7$ & $7.00\,10^{-2} \its 7$ \\
         72422 & $3.26\,10^{0} \its{13}$ &  $3.20\,10^{0} \its{10}$ & $3.21\,10^{0} \its{8}$ & $3.14\,10^{0} \its{8}$ \\
         648444 & $3.44\,10^2 \its{24}$ & $3.47\,10^2 \its{16}$ & $3.36\,10^2 \its 9$ & $3.31\,10^2 \its 8$ \\ \bottomrule
        \end{tabular}
        \caption{Exact-AMS}
    \end{subtable}
    \begin{subtable}{0.6\textwidth}
        \begin{tabular}{r|r|r|r|r}
            \toprule Dofs & \multicolumn{4}{c}{Solution time (linear iterations)} \\
         \midrule & \makecell[c]{$\Delta t=0.1$} & \makecell[c]{$\Delta t=0.05$} & \makecell[c]{$\Delta t=0.01$} & \makecell[c]{$\Delta t=0.005$} \\ 
         \cmidrule(lr){2-5} 349 &  $8.40\,10^{-3} \its{19}$  & $8.36\,10^{-3} \its{18}$ & $8.09\,10^{-3} \its{15}$ & $8.01\,10^{-3} \its{14}$ \\
         1404 &  $1.71\,10^{-2} \its{25}$ & $1.71\,10^{-2} \its{24}$ & $1.58\,10^{-2} \its{20}$ & $1.52\,10^{-2} \its{19}$ \\
         9156 & $8.54\,10^{-2}\its {33}$ & $7.98\,10^{-2} \its{30}$ & $7.22\,10^{-2} \its{26}$ & $6.83\,10^{-2} \its{24}$ \\
         72422 & $9.17\,10^{-1} \its{43}$ &  $7.74\,10^{-1} \its{36}$ & $6.46\,10^{-1} \its{28}$ & $6.22\,10^{-1} \its{27}$ \\
         648444 & $2.29\,10^1 \its{67}$ & $1.53\,10^1 \its{44}$ & $1.09\,10^1 \its{30}$ & $1.03\,10^1 \its{28}$ \\ \bottomrule
        \end{tabular}
        \caption{Jacobi-AMS}
    \end{subtable}
    \caption{Optimality test: Tetrahedral elements with VEM.}
    \label{tab:optimality-tetgen}
\end{table}
\begin{table}
    \centering
    \small
    \begin{subtable}{0.6\textwidth}
        \begin{tabular}{r|r|r|r|r}
            \toprule Dofs & \multicolumn{4}{c}{Solution time (linear iterations)} \\
        \midrule & \makecell[c]{$\Delta t=0.1$} & \makecell[c]{$\Delta t=0.05$} & \makecell[c]{$\Delta t=0.01$} & \makecell[c]{$\Delta t=0.005$} \\  
         \cmidrule(lr){2-5}  90 & $6.54\,10^{-3} \its 2$  & $6.69\,10^{-3} \its 2$ & $6.57\,10^{-3} \its 2$ & $6.76\,10^{-3} \its 2$ \\
         540 & $7.65\,10^{-3} \its 2$ & $7.78\,10^{-3} \its 2$ & $7.82\,10^{-2} \its 2 $ & $7.80\,10^{-3} \its 2$ \\
         3672 & $3.83\,10^{-2} \its 2 $ & $3.70\,10^{-2} \its 2 $ & $3.87\,10^{-2} \its 2 $ & $3.70\,10^{-2} \its 2 $ \\
         26928 & $2.48\,10^0 \its 2 $ &  $2.40\,10^0 \its 2 $ & $2.50\,10^0 \its 2 $ & $2.50\,10^0 \its 2 $ \\
         205920 & $2.31\,10^2 \its 2$ & $2.32\,10^2 \its 2$ & $2.32\,10^2 \its 2$ & $2.31\,10^2 \its 2$ \\ \bottomrule
        \end{tabular}
        \caption{Exact-exact}
    \end{subtable}
    \begin{subtable}{0.6\textwidth}
        \begin{tabular}{r|r|r|r|r}
            \toprule Dofs & \multicolumn{4}{c}{Solution time (linear iterations)} \\
          \midrule & \makecell[c]{$\Delta t=0.1$} & \makecell[c]{$\Delta t=0.05$} & \makecell[c]{$\Delta t=0.01$} & \makecell[c]{$\Delta t=0.005$} \\ 
         \cmidrule(lr){2-5} 90 & $1.15\,10^{-2}\its 2$  & $7.21\,10^{-3} \its 2 $ & $7.38\,10^{-3} \its 2 $ & $6.58\,10^{-3} \its 2 $ \\
         540 & $7.74\,10^{-3} \its 7 $ & $7.96\,10^{-3} \its 6 $ & $7.66\,10^{-2} \its 4 $ & $7.65\,10^{-3} \its 4$ \\
         3672 &  $3.24\,10^{-2} \its{12}$ & $3.19\,10^{-2} \its{10}$ & $3.02\,10^{-2} \its 6$ & $3.01\,10^{-2} \its 6$ \\
         26928 & $1.98\,10^0 \its{16}$ &  $1.95\,10^0 \its{13}$ & $1.92\,10^0 \its 9$ & $1.90\,10^0 \its 8$ \\
         205920 & $1.81\,10^2 \its{19}$ & $1.81\,10^2 \its{16}$ & $1.80\,10^2 \its{10}$ & $1.79\,10^2 \its 9$ \\ \bottomrule
        \end{tabular}
        \caption{Jacobi-exact}
    \end{subtable}   
    \begin{subtable}{0.6\textwidth}
        \begin{tabular}{r|r|r|r|r}
            \toprule Dofs & \multicolumn{4}{c}{Solution time (linear iterations)} \\
         \midrule & \makecell[c]{$\Delta t=0.1$} & \makecell[c]{$\Delta t=0.05$} & \makecell[c]{$\Delta t=0.01$} & \makecell[c]{$\Delta t=0.005$} \\ 
         \cmidrule(lr){2-5} 90 & $1.21\,10^{-2} \its 2$  & $7.39\,10^{-3} \its 2 $ & $7.79\,10^{-3} \its 2 $ & $1.21\,10^{-2} \its 2$ \\
         540 & $9.78\,10^{-3} \its 3$ & $1.00\,10^{-2} \its 3$ & $1.01\,10^{-2} \its 3$ & $1.27\,10^{-2} \its 3$ \\
         3672 & $4.05\,10^{-2} \its 5 $ & $3.84\,10^{-2} \its 4$ & $3.78\,10^{-2} \its 4$ & $3.82\,10^{-2} \its 4$ \\
         26928 & $8.79\,10^{-1} \its 8$ &  $8.27\,10^{-1} \its 6$ & $7.98\,10^{-1} \its 5$ & $7.70\,10^{-1} \its 5$ \\
         205920 & $6.09\,10^1 \its{14}$ & $5.89\,10^1 \its 9$ & $5.70\,10^1 \its 5$ & $5.72\,10^1 \its 5$ \\ \bottomrule
        \end{tabular}
        \caption{Exact-AMS}
    \end{subtable}
    \begin{subtable}{0.6\textwidth}
        \begin{tabular}{r|r|r|r|r}
            \toprule Dofs & \multicolumn{4}{c}{Solution time (linear iterations)} \\
         \midrule & \makecell[c]{$\Delta t=0.1$} & \makecell[c]{$\Delta t=0.05$} & \makecell[c]{$\Delta t=0.01$} & \makecell[c]{$\Delta t=0.005$} \\ 
         \cmidrule(lr){2-5} 90 &$7.73\,10^{-3} \its 2$  & $7.69\,10^{-3} \its 2$ & $7.39\,10^{-3} \its 2$ & $7.73\,10^{-3} \its 2$ \\
         540 & $1.05\,10^{-2} \its 7$ & $1.05\,10^{-2} \its 6$ & $1.02\,10^{-2} \its 4$ & $1.37\,10^{-2} \its 4$ \\
         3672 & $4.73\,10^{-2} \its{13}$ & $4.06\,10^{-2} \its{10}$ & $3.44\,10^{-2} \its 7$ & $3.36\,10^{-2} \its 7$ \\
         26928 & $4.59\,10^{-1} \its{20}$ &  $3.52\,10^{-1} \its{14}$ & $2.56\,10^{-1} \its 9$ & $2.60\,10^{-1} \its 9 $ \\
         205920 & $6.71\,10^0 \its{32}$ & $4.68\,10^0 \its{21}$ & $2.71\,10^0 \its{11}$ & $2.54\,10^0 \its{10}$ \\ \bottomrule
        \end{tabular}
        \caption{Jacobi-AMS}
    \end{subtable}
    \caption{Optimality test: Structured cubes mesh with VEM.}
    \label{tab:optimality-struct}
\end{table} 
\begin{table}
    \centering
    \small
    \begin{subtable}{0.6\textwidth}
        \begin{tabular}{r|r|r|r|r}
            \toprule Dofs & \multicolumn{4}{c}{Solution time (linear iterations)} \\
         \midrule & \makecell[c]{$\Delta t=0.1$} & \makecell[c]{$\Delta t=0.05$} & \makecell[c]{$\Delta t=0.01$} & \makecell[c]{$\Delta t=0.005$} \\  
         \cmidrule(lr){2-5} 110 & $7.04\,10^{-3} \its 2$  & $6.83\,10^{-3} \its 2$ & $6.99\,10^{-3} \its 2$ & $6.76\,10^{-3} \its 2$ \\
         700 & $1.29\,10^{-2} \its 2$ & $9.62\,10^{-3} \its 2$ & $9.62\,10^{-3} \its 2$ & $9.35\,10^{-3} \its 2$ \\
         4952 & $9.75\,10^{-2} \its 2$ & $9.68\,10^{-2} \its 2$ & $9.71\,10^{-2} \its 2$ & $9.71\,10^{-2} \its 2$ \\
         37168 & $6.68\,10^0 \its 2$ &  $6.69\,10^0 \its 2$ & $6.73\,10^0 \its 2$ & $6.59\,10^0 \its 2$ \\
         237690 & $3.96\,10^2 \its 2$ & $3.92\,10^2 \its 2$ & $3.97\,10^2 \its 2$ & $3.81\,10^2 \its 2$ \\ \bottomrule
        \end{tabular}
        \caption{Exact-exact}
    \end{subtable}
    \begin{subtable}{0.6\textwidth}
        \begin{tabular}{r|r|r|r|r}
            \toprule Dofs & \multicolumn{4}{c}{Solution time (linear iterations)} \\
         \midrule & \makecell[c]{$\Delta t=0.1$} & \makecell[c]{$\Delta t=0.05$} & \makecell[c]{$\Delta t=0.01$} & \makecell[c]{$\Delta t=0.005$} \\ 
         \cmidrule(lr){2-5} 110 & $7.03\,10^{-3} \its 8$  & $7.27\,10^{-3} \its 7 $ & $6.89\,10^{-3} \its 6$ & $6.76\,10^{-3} \its 6$ \\
         700 & $9.12\,10^{-3} \its{12}$ & $9.46\,10^{-3} \its{11}$ & $9.01\,10^{-3} \its 9$ & $9.04\,10^{-3} \its 8$ \\
         4952 & $9.44\,10^{-2} \its{16}$ & $9.36\,10^{-2} \its{14}$ & $9.04\,10^{-2} \its{11}$ & $8.95\,10^{-2} \its{10}$ \\
         37168 & $6.06\,10^0 \its{21}$ &  $5.97\,10^0 \its{16}$ & $5.91\,10^0 \its{12}$ & $5.78\,10^0 \its{11}$ \\
         237690 & $3.39\,10^2 \its{25}$ & $3.38\,10^2 \its{21}$ & $3.30\,10^2 \its{13}$ & $3.28\,10^2 \its{11}$ \\ \bottomrule
        \end{tabular}
        \caption{Jacobi-exact}
    \end{subtable}   
    \begin{subtable}{0.6\textwidth}
        \begin{tabular}{r|r|r|r|r}
            \toprule Dofs & \multicolumn{4}{c}{Solution time (linear iterations)} \\
          \midrule & \makecell[c]{$\Delta t=0.1$} & \makecell[c]{$\Delta t=0.05$} & \makecell[c]{$\Delta t=0.01$} & \makecell[c]{$\Delta t=0.005$} \\ 
         \cmidrule(lr){2-5} 110 & $7.69\,10^{-3} \its 2$  & $7.78\,10^{-3}\its 2$ & $7.85\,10^{-3}\its 2$ & $7.59\,10^{-3}\its 2$ \\
         700 & $1.17\,10^{-2} \its 3$ & $1.13\,10^{-2} \its 3$ & $1.11\,10^{-2} \its 3$ & $1.09\,10^{-2} \its 3$ \\
         4952 & $6.54\,10^{-2} \its 6$ & $6.02\,10^{-2} \its 5$ & $5.65\,10^{-2} \its 4$ & $5.59\,10^{-2} \its 4$ \\
         37168 & $1.69\,10^0 \its{11}$ &  $1.54\,10^0 \its 8$ & $1.40\,10^0 \its 5$ & $1.40\,10^0 \its 5$ \\
         237690 & $7.08\,10^1 \its{18}$ & $6.76\,10^1 \its{11}$ & $6.55\,10^1\its 6$ & $6.49\,10^1 \its 6$ \\ \bottomrule
        \end{tabular}
        \caption{Exact-AMS}
    \end{subtable}
    \begin{subtable}{0.6\textwidth}
        \begin{tabular}{r|r|r|r|r}
            \toprule Dofs & \multicolumn{4}{c}{Solution time (linear iterations)} \\
         \midrule & \makecell[c]{$\Delta t=0.1$} & \makecell[c]{$\Delta t=0.05$} & \makecell[c]{$\Delta t=0.01$} & \makecell[c]{$\Delta t=0.005$} \\ 
         \cmidrule(lr){2-5} 110 & $7.59\,10^{-3} \its 8$  & $7.91\,10^{-3} \its 7$ & $7.76\,10^{-3} \its 6$ & $7.55\,10^{-3} \its 6$ \\
         700 & $1.61\,10^{-2} \its{12}$ & $1.30\,10^{-2} \its{11} $ & $1.21\,10^{-2} \its 9$ & $1.28\,10^{-2} \its 8$ \\
         4952 & $8.61\,10^{-2} \its{16}$ & $7.76\,10^{-2} \its{14}$ & $6.56\,10^{-2} \its{11}$ & $6.18\,10^{-2} \its{10}$ \\
         37168 & $1.06\,10^0 \its{24}$ &  $7.86\,10^{-1} \its{17}$ & $5.97\,10^{-1} \its{12}$ & $5.49\,10^{-1} \its{11}$ \\
         237690 & $1.23\,10^1 \its{40}$ & $7.99\,10^0 \its{25}$ & $4.61\,10^0 \its{13}$ & $4.19\,10^0 \its{12}$ \\ \bottomrule
        \end{tabular}
        \caption{Jacobi-AMS}
    \end{subtable}
    \caption{Optimality test: Octahedral elements mesh with VEM.}
    \label{tab:optimality-octa}
\end{table} 
\begin{table}
    \centering
    \small
    \begin{subtable}{0.6\textwidth}
        \begin{tabular}{r|r|r|r|r}
            \toprule Dofs & \multicolumn{4}{c}{Solution time (linear iterations)} \\
        \midrule & \makecell[c]{$\Delta t=0.1$} & \makecell[c]{$\Delta t=0.05$} & \makecell[c]{$\Delta t=0.01$} & \makecell[c]{$\Delta t=0.005$} \\ 
         \cmidrule(lr){2-5} 206 & $7.11\,10^{-3} \its 2$  & $7.23\,10^{-3} \its 2$ & $7.28\,10^{-3} \its 2$ & $7.16\,10^{-3} \its 2$ \\
         1408 & $1.91\,10^{-2} \its 2$ & $1.88\,10^{-2} \its 2$ & $1.92\,10^{-2} \its 2$ & $1.88\,10^{-2} \its 2$ \\
         10376 & $8.46\,10^{-1} \its 2$ & $8.64\,10^{-1} \its 2$ & $8.66\,10^{-1} \its 2$ & $8.75\,10^{-1} \its 2$ \\
         79600 & $9.00\,10^1 \its 2$ & $8.98\,10^1 \its 2$ & $9.07\,10^1 \its 2$ & $8.85\,10^1 \its 2$ \\ \bottomrule
        \end{tabular}
        \caption{Exact-exact}
    \end{subtable}
    \begin{subtable}{0.6\textwidth}
        \begin{tabular}{r|r|r|r|r}
            \toprule Dofs & \multicolumn{4}{c}{Solution time (linear iterations)} \\
         \midrule & \makecell[c]{$\Delta t=0.1$} & \makecell[c]{$\Delta t=0.05$} & \makecell[c]{$\Delta t=0.01$} & \makecell[c]{$\Delta t=0.005$} \\ 
         \cmidrule(lr){2-5} 206 & $7.15\,10^{-3} \its{13}$  & $7.39\,10^{-3} \its{12}$ & $7.18\,10^{-3} \its{10}$ & $7.24\,10^{-3} \its 9$ \\
         1408 & $1.96\,10^{-2} \its{17}$ & $1.94\,10^{-2} \its{15}$ & $1.87\,10^{-2} \its{12}$ & $1.83\,10^{-2} \its{11}$ \\
         10376 & $9.69\,10^{-1} \its{23}$ & $8.99\,10^{-1} \its{19}$ & $8.60\,10^{-1} \its{14}$ & $8.90\,10^{-1} \its{13}$ \\
         79600 & $9.00\,10^1 \its{51}$ & $8.98\,10^1 \its{45}$ & $9.07\,10^1 \its{32}$ & $8.85\,10^1 \its{27}$ \\ \bottomrule
        \end{tabular}
        \caption{Jacobi-exact}
    \end{subtable}   
    \begin{subtable}{0.6\textwidth}
        \begin{tabular}{r|r|r|r|r}
            \toprule Dofs & \multicolumn{4}{c}{Solution time (linear iterations)} \\
          \midrule & \makecell[c]{$\Delta t=0.1$} & \makecell[c]{$\Delta t=0.05$} & \makecell[c]{$\Delta t=0.01$} & \makecell[c]{$\Delta t=0.005$} \\ 
         \cmidrule(lr){2-5} 206 & $8.70\,10^{-3} \its 3$  & $8.51\,10^{-3} \its 3$ & $8.47\,10^{-3} \its 2$ & $8.79\,10^{-3} \its 2$ \\
         1408 & $2.39\,10^{-2} \its 4$ & $2.32\,10^{-2} \its 4$ & $2.20\,10^{-2} \its 3$ & $2.15\,10^{-2} \its 3$ \\
         10376 & $2.92\,10^{-1} \its{10}$ & $2.42\,10^{-1} \its 7$ & $2.09\,10^{-1} \its 5$ & $1.96\,10^{-1} \its 4$ \\
         79600 & $8.62\,10^0 \its{17}$ & $7.54\,10^0 \its{11}$ & $6.67\,10^0 \its 6$ & $6.34\,10^0 \its 5$ \\ \bottomrule
        \end{tabular}
        \caption{Exact-AMS}
    \end{subtable}
    \begin{subtable}{0.6\textwidth}
        \begin{tabular}{r|r|r|r|r}
            \toprule Dofs & \multicolumn{4}{c}{Solution time (linear iterations)} \\
         \midrule & \makecell[c]{$\Delta t=0.1$} & \makecell[c]{$\Delta t=0.05$} & \makecell[c]{$\Delta t=0.01$} & \makecell[c]{$\Delta t=0.005$} \\ 
         \cmidrule(lr){2-5} 206 & $9.56\,10^{-3} \its{13}$  & $9.44\,10^{-3} \its{12}$ & $9.01\,10^{-3} \its{10}$ & $9.12\,10^{-3} \its 9$ \\
         1408 & $4.52\,10^{-2} \its{17}$ & $4.03\,10^{-2} \its{15}$ & $3.49\,10^{-2} \its{12}$ & $3.28\,10^{-2} \its{11}$ \\
         10376 & $4.11\,10^{-1} \its{24}$ & $3.49\,10^{-1} \its{20}$ & $2.61\,10^{-1} \its{15}$ & $2.46\,10^{-1} \its{13}$ \\
         79600 & $8.98\,10^0 \its{57}$ & $7.57\,10^0 \its{48}$ & $5.11\,10^0 \its{32}$ & $4.37\,10^0 \its{28}$ \\ \bottomrule
        \end{tabular}
        \caption{Jacobi-AMS}
    \end{subtable}
    \caption{Optimality test: Nonahedral elements mesh with VEM.}
    \label{tab:optimality-novedri}
\end{table} 
\begin{table}
    \centering
    \small
    \begin{subtable}{0.6\textwidth}
        \begin{tabular}{r|r|r|r|r}
            \toprule Dofs & \multicolumn{4}{c}{Solution time (linear iterations)} \\
        \midrule & \makecell[c]{$\Delta t=0.1$} & \makecell[c]{$\Delta t=0.05$} & \makecell[c]{$\Delta t=0.01$} & \makecell[c]{$\Delta t=0.005$} \\  
         \cmidrule(lr){2-5}  452 & $1.03\,10^{-2} \its 2$  & $8.70\,10^{-3} \its 2$ & $8.78\,10^{-3} \its 2$ & $8.64\,10^{-3} \its 2$ \\
         2047 & $4.56\,10^{-2} \its 2$ & $4.62\,10^{-2}\its 2$ & $4.72\,10^{-2}\its 2$ & $4.69\,10^{-2} \its 2$ \\
         17511 & $5.55\,10^{0} \its 2$ & $5.50\,10^{0} \its 2$ & $5.55\,10^{0} \its 2$ & $5.40\,10^{0} \its 2$ \\
         144862 & $6.12\,10^2 \its 2$ &  $6.04\,10^2 \its 2$ & $6.09\,10^2\its 2$ & $5.89\,10^2 \its 2$ \\
         292044 & $2.50\,10^3 \its 4$ & $2.51\,10^3 \its 2$ & $2.50\,10^3 \its 2$ & $2.41\,10^3 \its 2$ \\ \bottomrule
        \end{tabular}
        \caption{Exact-exact}
    \end{subtable}
    \begin{subtable}{0.6\textwidth}
        \begin{tabular}{r|r|r|r|r}
            \toprule Dofs & \multicolumn{4}{c}{Solution time (linear iterations)} \\
          \midrule & \makecell[c]{$\Delta t=0.1$} & \makecell[c]{$\Delta t=0.05$} & \makecell[c]{$\Delta t=0.01$} & \makecell[c]{$\Delta t=0.005$} \\ 
         \cmidrule(lr){2-5} 452 & $9.31\,10^{-3} \its{25}$  & $9.10\,10^{-3} \its{23}$ & $9.20\,10^{-3} \its{20}$ & $8.94\,10^{-3} \its{19}$ \\
         2047 & $5.01\,10^{-2} \its{20}$ & $4.88\,10^{-2} \its{19}$ & $4.77\,10^{-2} \its{15}$ & $4.82\,10^{-2} \its{14}$ \\
         17511 & $5.43\,10^{0} \its{27}$ & $5.46\,10^{0} \its{24}$ & $5.37\,10^{0} \its{19}$ & $5.25\,10^{0} \its{17}$ \\
         144862 & $5.85\,10^2 \its{34}$ &  $5.76\,10^2 \its{29}$ & $5.71\,10^2 \its{21}$ & $5.58\,10^2 \its{19}$ \\
         292044 & $2.40\,10^3 \its{39}$ & $2.39\,10^3 \its{32}$ & $2.39\,10^3 \its{22}$ & $2.28\,10^3 \its{21}$ \\ \bottomrule
        \end{tabular}
        \caption{Jacobi-exact}
    \end{subtable}   
    \begin{subtable}{0.6\textwidth}
        \begin{tabular}{r|r|r|r|r}
            \toprule Dofs & \multicolumn{4}{c}{Solution time (linear iterations)} \\
          \midrule & \makecell[c]{$\Delta t=0.1$} & \makecell[c]{$\Delta t=0.05$} & \makecell[c]{$\Delta t=0.01$} & \makecell[c]{$\Delta t=0.005$} \\ 
         \cmidrule(lr){2-5} 452 & $4.93\,10^{-2} \its{57}$  & $4.09\,10^{-2} \its{43}$ & $2.47\,10^{-2} \its{20}$ & $2.08\,10^{-2} \its{15}$ \\
         2047 & $2.59\,10^{-1} \its{57}$ & $1.68\,10^{-1} \its{35}$ & $8.02\,10^{-2} \its{13}$ & $6.36\,10^{-2} \its 9$ \\
         17511 & -- & $3.76\,10^{1} \its{614}$ & $1.25\,10^{1} \its{208}$ & $8.38\,10^{0} \its{142}$ \\
         144862 & -- &  -- & -- & $5.66\,10^2 \its{824}$ \\
         292044 & -- & -- & -- & -- \\ \bottomrule
        \end{tabular}
        \caption{Exact-AMS}
    \end{subtable}
    \begin{subtable}{0.6\textwidth}
        \begin{tabular}{r|r|r|r|r}
            \toprule Dofs & \multicolumn{4}{c}{Solution time (linear iterations)} \\
         \midrule & \makecell[c]{$\Delta t=0.1$} & \makecell[c]{$\Delta t=0.05$} & \makecell[c]{$\Delta t=0.01$} & \makecell[c]{$\Delta t=0.005$} \\ 
         \cmidrule(lr){2-5} 452 & $9.47\,10^{-2} \its{127}$  & $6.38\,10^{-2} \its{82}$ & $3.63\,10^{-2} \its{37}$ & $3.14\,10^{-2} \its{30}$ \\
         2047 & $4.06\,10^{-1} \its{96}$  & $2.61\,10^{-1} \its{60}$ & $1.32\,10^{-1} \its{28}$ & $1.13\,10^{-1} \its{23}$ \\
         17511 & -- & -- & $1.91\,10^{1} \its{340}$ & $1.32\,10^{1} \its{242}$ \\
         144862 & -- & -- & -- & -- \\
         292044 & -- & -- & -- & -- \\ \bottomrule
        \end{tabular}
        \caption{Jacobi-AMS}
    \end{subtable}
    \caption{Optimality test: Irregular polyhedral elements mesh with VEM.}
    \label{tab:optimality-voro}
\end{table} 

\subsection{Scalability tests}
In applications, the resulting meshes will be fairly regular, thus we perform the strong scalability tests on the tetrahdral, structured, octahedral and nonahedral meshes for a fixed timestep of $\tau=0.05$, as it was observed in the optimality tests to be the largest time step to yield a robust preconditioner. The results are shown in Table \ref{tab:scalability-vem}. We note that the results are in general positive, as all methods present only a very mild deterioration in performance as the number of CPU cores increases, and the decrease in solution time correspongs to the increase of CPUs up to 8 cores, where then communication times hinder the scalability of the preconditioner. If we denote by $T_p$ the solution time required by a method using $p$ processors, we can define the speed-up and efficiency as $S_p=T_1/T_p$ and $E_p=T_1/pT_p$ respectively, which we plot in Figre \ref{fig:speedup-eff}. In green we show the ideal curves, where indeed we see a deterioration of the performance as we increase the number of CPU cores. The effect on the tetrahedral mesh is not as strong on the tetrahedral mesh, this because we considered a much bigger problem for that case (roughly 650k DoFs instead of 240k-300k for the others).

\begin{table}
    \centering
    \small
    \begin{tabular}{r|r|r|r|r|r|r|r}
        \toprule & \multicolumn{6}{c}{Solution time (linear iterations)} \\
         \midrule & DoFs  & \makecell[c]{1 CPU} & \makecell[c]{2 CPU} & \makecell[c]{4 CPU} & \makecell[c]{8 CPU} & \makecell[c]{16 CPU} & \makecell[c]{32 CPU}  \\
         \cmidrule(lr){2-8} Tet & 648444 & $48.07 \its{44}$ & $21.61 \its{45}$ & $13.46 \its{45}$ & $6.70 \its{45}$ & $3.72 \its{45}$ & $ 2.22 
         \its{45}$ \\
         Struct & 291708 & $23.27\its{22}$ & $13.35 \its{24}$ & $6.97 \its{25}$ & $4.06 \its{27}$ & $2.19 \its{27}$ & $ 1.51 \its{28}$ \\
         Octa & 237690 & $22.95 \its{25}$ & $12.78 \its{26}$ & $6.81 \its{27}$ & $3.83 \its{28}$ & $2.15 \its{28}$ & $1.47 \its{28}$ \\ 
         Nona & 264888 & $37.60 \its{34}$ & $20.95 \its{36}$ & $13.26 \its{38}$ & $7.18 \its{39}$ & $3.34 \its{39}$ & $2.57 \its{39}$\\ \bottomrule
    \end{tabular}
    \caption{Scalability test: VEM.}
    \label{tab:scalability-vem}
\end{table}
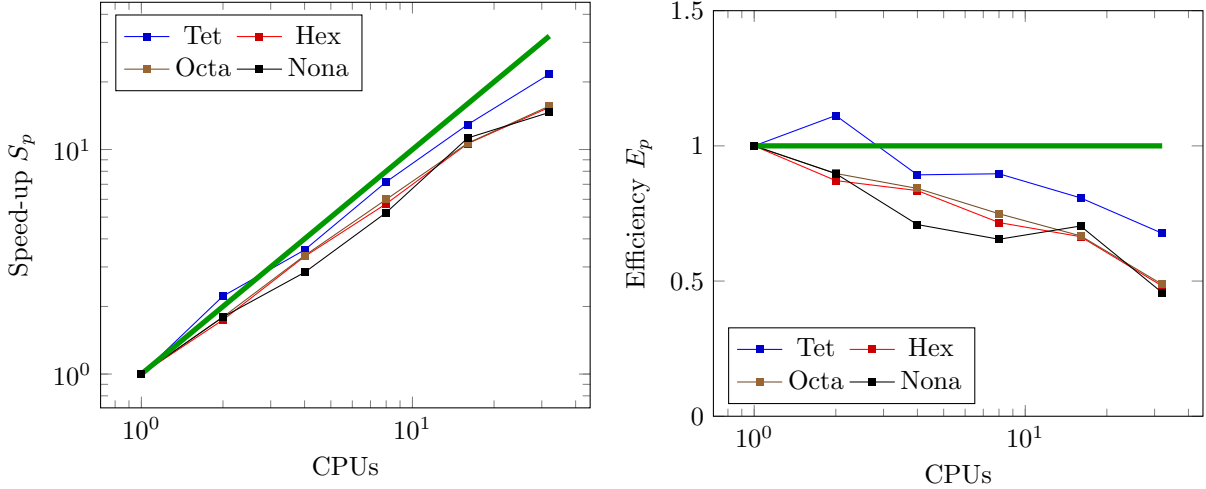
\begin{figure}
    \pgfplotsset{every axis plot/.append style={mark options={draw opacity=0.0, fill opacity=1.0, mark size=1.5pt, mark=square*, }} }
    \centering
    \begin{subfigure}{0.49\textwidth}
        \begin{tikzpicture}
    \begin{loglogaxis}
        [width=\textwidth,
        xlabel=CPUs,
        ylabel=Speed-up $S_p$,
        legend columns=2, 
        legend pos=north west,
        ]

        \addplot  coordinates { 
            (1,  1) 
            (2,  48.07 / 21.61) 
            (4,  48.07 / 13.46) 
            (8,  48.07 / 6.70) 
            (16, 48.07 / 3.72) 
            (32, 48.07 / 2.22) 
            };
        \addplot  coordinates { 
            (1,  1) 
            (2,  23.27 / 13.35) 
            (4,  23.27 / 6.97) 
            (8,  23.27 / 4.06) 
            (16, 23.27 / 2.19) 
            (32, 23.27 / 1.51) 
            };
        \addplot  coordinates { 
            (1,  1) 
            (2,  22.95 / 12.78) 
            (4,  22.95 / 6.81) 
            (8,  22.95 / 3.83) 
            (16, 22.95 / 2.15) 
            (32, 22.95 / 1.47) 
            };
        \addplot  coordinates { 
            (1,  1) 
            (2,  37.60 / 20.95) 
            (4,  37.60 / 13.26) 
            (8,  37.60 / 7.18) 
            (16, 37.60 / 3.34) 
            (32, 37.60 / 2.57) 
            };
        \addplot[line width=2pt, color=green!60!black] coordinates{
            (1,1)
            (2,2)
            (4,4)
            (8,8)
            (16,16)
            (32,32)};
        
    \legend{Tet, Hex, Octa, Nona}
    \end{loglogaxis}
    \end{tikzpicture}
     
    \end{subfigure}
    \begin{subfigure}{0.49\textwidth}
        \begin{tikzpicture}
    \begin{semilogxaxis}
        [width=\textwidth,
        xlabel=CPUs,
        ymin=0.,
        ymax=1.5,
        ylabel=Efficiency $E_p$,
         legend columns=2, 
         legend pos=south west
        ]

        \addplot  coordinates { 
            (1,  1) 
            (2,  48.07 / 2  / 21.61) 
            (4,  48.07 / 4  / 13.46) 
            (8,  48.07 / 8  / 6.70) 
            (16, 48.07 / 16 / 3.72) 
            (32, 48.07 / 32 / 2.22) 
            };
        \addplot  coordinates { 
            (1,  1) 
            (2,  23.27 / 2  /13.35) 
            (4,  23.27 / 4  /6.97) 
            (8,  23.27 / 8  /4.06) 
            (16, 23.27 / 16 /2.19) 
            (32, 23.27 / 32 /1.51) 
            };
        \addplot  coordinates { 
            (1,  1) 
            (2,  22.95 / 2/12.78)
            (4,  22.95 / 4/6.81)
            (8,  22.95 / 8/3.83)
            (16, 22.95 / 16/2.15)
            (32, 22.95 / 32/1.47)
            };
        \addplot  coordinates { 
            (1,  1) 
            (2,  37.60 / 2/20.95)
            (4,  37.60 / 4/13.26)
            (8,  37.60 / 8/7.18)
            (16, 37.60 / 16/3.34)
            (32, 37.60 / 32/2.57)
            };
        \addplot[line width=2pt, color=green!60!black] coordinates{
            (1,1)
            (2,1)
            (4,1)
            (8,1)
            (16,1)
            (32,1)};
    \legend{Tet, Hex, Octa, Nona}
    \end{semilogxaxis}
    \end{tikzpicture}
     
    \end{subfigure}
    \caption{Scalability tests, on the left we show the speed-up $S_p$, on the right the parallel efficiency $E_p$. Ideal curves (linear speed-up, perfect efficiency) are depcited in green to better appreciate the deterioration of performance when increasing the number of processors.}
    \label{fig:speedup-eff}
\end{figure}

\subsection{Sensitivity on stabilization in Schur operator}
As shown in Section \ref{sec:vem}, VEM has an intrinsic stabilization block required for the invertibility of the discrete operator. Still, the magnitude of the stabilization can be modulated through a parameter, which can have a great impact on the conditioning of the resulting problem, as studied in \cite{antonietti.2021}. The sparse representation of the Schur complement operator we are considering is not the same matrix as the one that appears in the problem formulation, so it is not clear how much stabilization is required for this block. This is of course true not only here, but it holds for any Schur complement preconditioner used in a problem with a block structure as this one. This motivated this preliminary investigation on the effect of the Schur stabilization on the performance of the preconditioner.

We show the results of this study in Table \ref{tab:stabilization}. We note that in all but the irregular case, the lowest iteration count is obtained by using the standard value $\alpha=1$. Still, this trend does not hold for the irregular case, where stabilization hinders performance. This last result is rather counter-intuitive, as a VEM block without stabilization is not invertible in theory. This is an interesting phenomenon that had not been observed before, and will indeed be the subject of future work. 
\begin{table}[ht]
    \centering
    \begin{tabular}{c|c||r|r|r|r|r|r}
        \toprule Mesh type & DoFs & $\alpha=0.001$ & $\alpha=0.01$ & $\alpha=0.1$ & $\alpha=1$ & $\alpha=10$ & $\alpha=100$\\
        \cmidrule(lr){1-2}\cmidrule(lr){3-8} Tet & 72422 & 36 & 36 & 36 & 36 & 36 & 36 \\ 
        Struct & 88200 & 33 & 32 & 28 & \textbf{20}& 40 & 109\\
        Octa & 37168 & 31 & 31 & 26 & \textbf{18} & 31 & 101\\
        Nona & 79600 & 63 & 59 & 41 & \textbf{28} & 55 & 151\\
        Voro & 17511 & \textbf{646} & 844 & 1018 & 1175 & -- & --\\\bottomrule
    \end{tabular}
    \caption{Linear iterations incurred when varying the stabilization parameter in the Schur complement operator. Tests were performed with a timestep of $\tau=0.05$. The lowest iteration case is highlighted with bold fonts, except for the tetrahedral mesh, where stabilization plays no role.}
    \label{tab:stabilization}
\end{table}
\section{Conclusions}\label{sec:conclusions}
This work presents an extension of an optimal FEM preconditioner to VEM, for a two field formulation of Maxwell's equations based on the Schur complement. Its performance is robust and scalable for a sufficiently (but not too much) small time step, with usually only a mild deterioration in the GMRES iterations. The source of deterioration for most cases is the dominance of the $\curl \curl$ operator in the Schur complement. When the meshes become too complex, this becomes the main driver of deterioration, as shown by the performance of the Jacobi preconditioner on the nonahedral mesh, or by the inability of the AMS preconditioner to function properly in the irregular mesh, which presents highly distorted elements. In all working cases, performance issues can be alleviated by reducing the time step, where the values required for an adequate performance are not too small, as we have observed a value of $0.05$ seconds to suffice. 

A comparison with finite elements, where the preconditioner used is known to be optimal, shows many interesting and still unresolved issues in the numerical approximation of PDEs using VEM. We highlight two: the first one is that the relationship between the performance and the elements used by the mesh are not trivial. Indeed, we expected the performance to deteriorate in the order tetrahedra, hexahedra, octahedra, nonahedra, and irregular, but the order between the first four is not like so. The second one is the preconditioning of the mass block, where usually a Jacobi preconditioner suffices, at least for first order elements. Instead, we observed that this approach was not effective with the nonahedral mesh. The latter point is fundamental for the simulation of time-dependent problems.

\section*{Acknowledgments}
This research was supported as follows: NAB, FD and SS were supported by INDAM GNCS.

\bibliographystyle{plain}
\bibliography{bibliography}

\end{document}